\RequirePackage[2019/01/01]{latexrelease}
\documentclass{autart}

\usepackage[ansinew]{inputenc}
\usepackage{graphicx}
\usepackage[T1]{fontenc}
\usepackage{amsmath}
\usepackage{amssymb}
\usepackage{bm}
\usepackage{bbold}
\usepackage{siunitx}
\usepackage{lipsum}
\usepackage{color}
\usepackage{cases}
\usepackage{IEEEtrantools}	
\usepackage{dsfont}
\usepackage{phonenumbers}

\newcommand{\mv}[1]{\boldsymbol{#1}}			
\newcommand{\mm}[1]{\mv{#1}} 		

\newcommand{\nullvec}[1]{\mathbb 0_{#1}}

\newcommand{\rank}{\text{rank}}

\newcommand{\dt}{\mathrm{d}t}

\newcommand{\R}{\mathds{R}}

\newcommand{\Rx}{\R^n}
\newcommand{\Ru}{\R^m}

\newcommand{\uss}{\mv{u}}						
\newcommand{\x}{\mv{x}}							
\newcommand{\Jphs}{\mm{J}}						
\newcommand{\Rphs}{\mm{R}}						
\newcommand{\dH}{\frac{\partial H}{\partial \x}}				
\newcommand{\dHtop}{\frac{\partial^\top H}{\partial \x^\top}}
\newcommand{\gphs}{\mm{G}}						

\newcommand{\tini}{t_0}							
\newcommand{\Ham}{\mathcal{H}}					
\newcommand{\argmin}{\arg\,\min}				
\newcommand{\Ropt}{\mm{S}(\bm x)}						
\newcommand{\Qopt}{r(\x)}						

\newcommand{\Zx}{\mathcal{X}_{\gphs}}								
\newcommand{\la}{\Upsilon}									
\newcommand{\clf}{\Vx}

\newcommand{\vtf}{\frac{\partial^\top V}{\partial \bm x^\top} \mv{f}(\bm x)}

\newcommand{\Sx}{\frac{\partial ^\top V}{\partial \bm x^\top} \gphs(\bm x) \mm{S} \gphs^\top(\bm x) \clf}

\newcommand{\fl}{{f}_\la}
\newcommand{\Ql}{{Q}_\la}
\newcommand{\Sl}{{S}_\la}
\newcommand{\Vx}{\frac{\partial V}{\partial \x}}   						   

\newcommand{\AnzahlGewichte}{r}
\newcommand{\Rw}{\R^{\AnzahlGewichte}}

\newcommand{\w}{\mv{w}}

\newcommand{\Qu}{\mathcal{Q}}

\definecolor{emerald}{rgb}{0,0.65,0.42}

\definecolor{blue}{rgb}{0,0,0.8}

\newcounter{subeq}

\begin{document}
\begin{frontmatter}
\title{Optimal Control of Port-Hamiltonian Systems: \\A Time-Continuous Learning Approach\thanksref{footnoteinfo}\thanksref{DFG}}
\thanks[footnoteinfo]{Corresponding author L.~Kölsch. Tel. \phonenumber[foreign]{0721 608}{43237}. 
	Fax \phonenumber[foreign]{0721 608}{42707}.}
\thanks[DFG]{This work was funded by the Deutsche Forschungsgemeinschaft (DFG, German Research Foundation)---project number 360464149.}
\author[irs]{Lukas Kölsch}\ead{lukas.koelsch@kit.edu},
\author[irs]{Pol Jané Soneira}\ead{pol.jane@kit.edu},
\author[irs]{Felix Strehle}\ead{felix.strehle@kit.edu},
\author[irs]{Sören Hohmann}\ead{soeren.hohmann@kit.edu}
\address[irs]{Institute of Control Systems, \\ Karlsruhe Institute of Technology (KIT),\\ Kaiserstr. 12, 76131 Karlsruhe,\\ Germany}
\begin{abstract}
Feedback controllers for port-Hamiltonian systems reveal an intrinsic inverse optimality property 	
since each passivating state feedback controller is optimal with respect to some specific performance index.
Due to the nonlinear port-Hamiltonian system structure, however, explicit (forward) methods for optimal control of port-Hamiltonian systems require the generally intractable analytical solution of the Hamilton-Jacobi-Bellman equation.
Adaptive dynamic programming methods provide a means to circumvent this issue. However, the few existing approaches for port-Hamiltonian systems hinge on very specific sub-classes of either performance indices or system dynamics or require the intransparent guessing of stabilizing initial weights. 
In this paper, we contribute towards closing this largely unexplored research area by proposing a time-continuous adaptive feedback controller for the optimal control of general time-continuous input-state-output port-Hamiltonian systems with respect to general Lagrangian performance indices.
Its control law implements an online learning procedure which uses the Hamiltonian of the system as an initial value function candidate.
The time-continuous learning of the value function is achieved by means of a certain Lagrange multiplier that allows to evaluate the optimality of the current solution.
In particular, constructive conditions for stabilizing initial weights are stated and
asymptotic stability of the closed-loop equilibrium is proven.
Our work is concluded by simulations for exemplary linear and nonlinear optimization problems which demonstrate asymptotic convergence of the controllers resulting from the proposed online adaptation procedure.
\end{abstract}
\end{frontmatter}
\section{Introduction}
In recent years, systematic modeling of dynamical multi-physics systems in port-Hamiltonian (pH) form has become increasingly popular in a wide range of applications such as acoustics \cite{Falaize.2016}, aerospace \cite{Aoues2019}, robotics \cite{Groothius,Macchelli2007,Macchelli2009}, power electronics \cite{Bergna2017,Cupelli2019}, and energy systems \cite{Fiaz.2013,StrehlePfeiferKoelsch2018}, just to name a few.
Due to their specific geometric, energy-based structure with power-conjugated port variables and a Hamiltonian representing the total stored energy in the system, pH systems constitute inherently passive systems \cite[Chapter~6-7]{vanderSchaft.2017}.
Thus, pH system representations are well suited for system analysis and control design based on passivity arguments. 
Besides feedback interconnection of passive systems as in Control by Interconnection (CbI), a common theme in standard passivity-based control is the passivation and asymptotic stabilization (under an additional detectability condition) via static state feedback $\bm u(\bm x)=-\bm k(\bm x)$ \cite[Chapter~2.4]{Sepulchre.1997},\cite[Chapters~5,7]{vanderSchaft.2017}, \cite{Ortega.2004,Ortega.2008}.

For the passivity-based state feedback controller design, an inverse optimality property can be characterized in both continuous \cite[p.~ 107ff.]{Sepulchre.1997},\cite[Theorem~3.5.1]{vanderSchaft.2017} and discrete time \cite{Monaco2015}.
For time-continuous, input-affine nonlinear systems
\begin{align}
\dot{\bm x}(t)= \bm f(\bm x(t)) + \bm g(\bm x(t)) \bm u(t), \label{input-affine-dynamics}
\end{align} 
it states that an input $\bm u(\bm x) =- \bm k(\bm x)$ is optimally stabilizing with respect to the specific performance index
\begin{align}
J(\bm x(t), \bm u(t)) &= \frac 12 \int_{t_0}^\infty \; \ell(\bm x(t)) + \bm u^\top(t) \bm u(t) \; \mathrm dt
\end{align}
with
\begin{align}
\ell(\bm x(t)) &= - \frac{\partial^\top V(\bm x(t))}{\partial \bm x} \bm f(\bm x(t)) + \frac 12 \bm k^\top (\bm x(t)) \bm k(\bm x(t))
\end{align}
if and only if the open-loop system is output feedback passive with positive-definite storage function $V(\bm x(t))$ and provided that some detectability condition is fulfilled, see \cite[p.~ 54ff.]{vanderSchaft.2017},\cite[Theorem~3.30]{Sepulchre.1997} for a detailed discussion.

However, classical, passivity-based control methods are in general not designed for optimal control problems in a more practical setup, where the ``forward'' solution to a given, arbitrary optimization problem with a general Lagrangian performance index 
is sought. 
In this paper, we address such general optimization problems subject to input-state-output port-Hamiltonian system (ISO-PHS) dynamics, this is
\begin{subequations}\label{original-OP}
	\begin{IEEEeqnarray}{C'l}\label{eq:guetemass-single-player}
		\min_{\uss(t)} &\frac{1}{2} \int_{\tini}^{\infty} \big(r(\bm x(t)) + (\uss(t))^\top \Ropt \uss(t) \big) \dt	\\
		\text{s.t.} &\dot{\x}(t) = (\Jphs(\bm x(t)) - \Rphs(\bm x(t)))\frac{\partial H(\bm x(t))}{\partial \bm x(t)} \nonumber \\ 
		& \quad \quad + \gphs(\bm x(t))\uss (t). \label{eq:systemdynamik_single_player}
	\end{IEEEeqnarray}
\end{subequations}
with state vector $\x \in \Rx$, input vector $\bm u \in \mathds R^m$, skew-symmetric interconnection matrix $\bm J\in \R^{n\times n} $, positive-semidefinite dissipation matrix $\bm R \in \R^{n\times n}$, input matrix $\bm G \in \mathds R^{n\times p}$, positive-definite Hamiltonian $H:\Rx \to \R_{\geq0}$, positive-definite $r: \Rx \to \R_{>0}$ and positive-definite $\bm S \in \R^{m\times m}$.
The system \eqref{eq:systemdynamik_single_player} is equipped with the passive output
\begin{align}
\bm y(t) = \bm G^\top(\bm x) \frac{\partial H(\bm x(t))}{\partial \bm x(t)}. \label{eq:output}
\end{align}
The following section provides an overview of current research in optimization-based control for ISO-PHSs.
\subsection{Related Work}


In the case of linear ISO-PHSs, the Hamiltonian $H(\bm x)$ is quadratic, which allows to calculate an optimal controller by using state-dependent Riccati equations \cite{Pei.2011}. In \cite{MarquisFavre.2008,Mouhib.2009,Gerelli.2009}, 
the necessary conditions which follow from Pontryagin's Maximum Principle are used to derive an explicit expression for the optimal feedback controller, provided the Hamiltonian of the system is quadratic. The authors in \cite{Wu.2018} provide full- and reduced-order LQR controllers for linear ISO-PHSs.
Further extensions of LQ-optimal control for pH systems are given for stochastic or infinite-dimensional spaces \cite{Lamoline.2018} and boundary control systems \cite{Liu.2020}.

While there is a rich theory available for the linear case, the general solution of optimal control problems for nonlinear ISO-PHSs remains challenging due to the necessity of explicitly solving the Hamilton-Jacobi-Bellman Equation (HJBE), which is a nonlinear PDE and thus hard to solve. 
This issue can be circumvented by applying adaptive dynamic programming (ADP) methods.
If the performance index of the optimal control problem has a specific structure and the system dynamics is given by a Hamiltonian system with controlled Hamiltonian $H(\bm x, \bm u)$, iterative learning control \cite{Fujimoto.2003,Fujimoto.2003b} and iterative feedback tuning methods \cite{Fujimoto.2008} have been proposed.
For the specific sub-class of fully actuated mechanical pH systems, the authors in \cite{Okura.2020} propose an adaptive path-following controller from a training trajectory using Bayesian estimation. 
\cite{Nageshrao.2014} and \cite{Sprangers.2015} use actor-critic reinforcement learning schemes to minimize the error between the resulting closed-loop system and a given desired closed-loop system without the need of explicitly solving the matching PDE of the employed passivity-based controller.
However, these approaches suffer from the dissipation obstacle, according to which the Hamiltonian can only be shaped for coordinates that are not affected by physical damping.
Thus, the Hamiltonian of the desired closed-loop system can not be freely chosen.
A profound overview on recent adaptive and learning-based control methods for pH systems can be found in \cite{Nageshrao.2016}.
If the optimization problem is convex and the performance index depends only on the final value of $\bm x$, applying the primal-dual gradient method results in a controller which is again port-Hamiltonian \cite{Stegink.2015}. 
This method is convenient for a wide range of practical applications and easy to implement \cite{Stegink.2017,Koelsch.2020}.
However, it does not allow to take the transient behaviour of the state or input trajectories into account in the optimization problem.

For the more general class of time-continuous input-affine nonlinear systems, ADP methods \cite{Vamvoudakis.2010,Jiang.2013,Bian.2016} are proposed where the optimal value function $V(\bm x)$ is iteratively found using a weighted sum of basis functions.
However, a proper set of initial weights leading to a stabilizing controller has to be found by educated guessing.

Overall, in the existing literature on optimal control of ISO-PHSs either the performance index or the system dynamics or both remain limited to some very specific sub-classes.
Likewise, ADP methods for ISO-PHSs as well as for the more general class of time-continuous, input-affine nonlinear systems are usually not constructive in the sense that they require the intransparent guessing of initial weights for a stabilizing value function candidate.
To the best of the authors' knowledge there exist no explicit control schemes for the dynamic optimal control problem \eqref{original-OP} with generalized Lagrangian performance index and a general ISO-PHS. 
\subsection{Main Contribution}
%
In this paper, we address this issue by developing a time-continuous adaptive feedback control strategy for the dynamic optimization problem \eqref{original-OP}.
The initial step of our design is based on a trick originally outlined in \cite{Sackmann2000}. By multiplying the system dynamic constraints \eqref{eq:systemdynamik_single_player} with the gradient of a control-Lyapunov function (CLF) $V(\bm x)$, we obtain a \emph{Modified Optimal Control} (MOC) problem that allows for an analytical solution of an asymptotically stabilizing $\uss(\x(t))=-\bm k(\bm x(t))$.
However, as a consequence of the modification, the MOC control law is optimal to an unintentionally modified objective function.

To achieve optimality with respect to the original optimization problem \eqref{original-OP}, we extend the MOC law by a gradient-based adaptation for the CLF. This ensures convergence of the CLF to the value function of \eqref{original-OP} and results in an explicit controller for optimization problem \eqref{original-OP}.
Furthermore, we derive necessary and sufficient conditions for Hamiltonians to be CLFs and show that if the Hamiltonian is a CLF, we are able to provide stabilizing initial weights for our adaptation strategy. Finally, we prove (asymptotic) stability of the closed-loop system equilibrium.
\subsection{Paper Organization}
The outline of this paper is as follows. 
In Section 2, we summarize the main results of \cite{Sackmann2000} and set them in the context of ISO-PHSs. 
After presenting an analytical solution of the MOC problem, we derive necessary and sufficient condition for the Hamiltonian $H(\bm x)$ being a CLF.
In Section 3, the modified optimal controller is enhanced by a learning procedure in order to achieve optimality with regard to the original problem \eqref{original-OP}. The resulting adaptive optimal controller is proven to be (asymptotically) stable.
Section 4 presents both linear and nonlinear examples showing the asymptotic convergence and optimality of the proposed control law.
A discussion and outlook on further research directions in Section 5 concludes our work.  
\subsection{Notation}
Both vectors and matrices are written in boldface.
All vectors defined in the paper are column vecotrs $\bm a = \mathrm{col}_i\{a_i\}=\mathrm{col}\{a_1,a_2,\ldots\}$ with elements $a_i$, $i=1,2,\ldots$.
All-zeros and all-ones vectors of dimension $n$ are denoted by $\nullvec{n}$ and $\mathds 1_n$, respectively.
The $(n \times n)$-identity matrix is denoted by $\bm I_n$.
Positive-semidefinite and -definite matrices or functions are denoted by $\succeq 0$ and $\succ 0$, respectively.
Equilibrium variables of the state $\bm x(t) \in \Rx$ are marked with a star
and shifted values with respect to an equilibrium are marked with a tilde, i.e. $\widetilde{\bm x}(t) = \bm x(t) - \bm x^\star$.
For vectors $\bm a, \bm b$ of the same size we write $\bm a \geq \bm b$ if each component in $\bm a$ is greater than or equal to the corresponding component in $\bm b$.
The set $\mathcal B(\bm x(t),\varepsilon)$ denotes a ball of radius $\varepsilon > 0$ around $\bm x(t)$.
To allow distinction from the Hamiltonian $H$ of the ISO-PHS, the Hamiltonian function of the optimization problem is denoted by $\mathcal H$.
For clarity of presentation, the time dependence $(t)$ of the variables is not explicitly mentioned anymore, unless it is essential for transparency of the statements.
\section{Modified Optimal Control for Port-Ha\-mil\-to\-nian Systems}
The starting point of our work is based on the MOC approach originally outlined in \cite{Sackmann2000}. If the constraints \eqref{eq:systemdynamik_single_player} are projected via the gradient of a CLF onto $\mathds R$, an MOC problem
\begin{subequations}\label{eq:morUrversion}
	\begin{IEEEeqnarray}{C,l}
		\min_{\uss} &\frac{1}{2} \int_{\tini}^{\infty} \big(r(\bm x) + (\uss)^\top \Ropt \uss \big) \dt	\\
		\mathrm{s.t.}  &\dot{V}(\x) = \frac{\partial^\top V(\bm x)}{\partial \bm x^\top} \Big(
		(\Jphs(\bm x) - \Rphs(\bm x))\frac{\partial H(\bm x)}{\partial \bm x} \nonumber \\ 
		& \quad \quad \quad + \gphs(\bm x)\uss \Big), 
		\label{eq:systemdynamik-MOD}
	\end{IEEEeqnarray}
\end{subequations}
arises, which can be solved analytically by a static state feedback $\bm u(\bm x) = - \bm k(\bm x)$.
However, the resulting feedback is in general not optimal with respect to the original problem \eqref{original-OP} and a suitable CLF must be found. While the former has not been addressed in literature yet, the latter is a stumbling block for general input-affine nonlinear systems \cite{Kokotovic.2001}. In the case of ISO-PHSs, however, the Hamiltonian $H(\bm x)$ presents a natural CLF candidate.

After revising the MOC approach \cite{Sackmann2000} in Section \ref{ch:introductionMOR}, we show in Section \ref{ch:PHSCLF}, to what extent CLFs can be found naturally in ISO-PHSs by using the Hamiltonian $H(\bm x)$. 
Finally, in Section \ref{ch:MOROC} we discuss the relationship between MOC and Optimal Control by deriving necessary and sufficient conditions on the Modified Optimal Controller to be an optimizer of the original problem \eqref{original-OP}.
This will form the basis for the following derivation of our learning procedure.
\subsection{Introduction to Modified Optimal Control}\label{ch:introductionMOR}
\begin{defn}[Control-Lyapunov function]\cite[p.~46]{Freeman.2008}\label{def:clf}
	A CLF for the system
	\begin{align}
		\dot{\x} = \mv{f}(\x, \uss), \quad \x \in \Rx, \quad \uss \in \Ru
	\end{align} 
	with $\bm f(\nullvec{n}, \nullvec{m})= \nullvec{n}$ is a radially unbounded, positive-definite function $V: \mathds R^n \to \mathds R$, fulfilling 
	\begin{align}\label{eq:bedingung_clf}
		\forall \, \x \neq \nullvec{n}: && \inf_{\uss}\bigg\{ \frac{\partial^\top V}{\partial \x^\top} \mv{f}(\x, \uss) \bigg\} < 0.
	\end{align}
\end{defn}
For input-affine nonlinear systems \eqref{input-affine-dynamics},
condition \eqref{eq:bedingung_clf}
is equivalent to \cite[p.~641]{Kokotovic.2001}
\begin{align}\label{def:clf:affine}
	\frac{\partial^\top V}{\partial \bm x^\top} \gphs(\x) = \nullvec{m} \quad \Rightarrow \quad \frac{\partial^\top V}{\partial \x^\top} \mv{f}(\bm x) \begin{cases}
	< 0, \quad  \x \neq \nullvec{n}, \\ = 0, \quad \x = \nullvec{n},
	\end{cases}
\end{align}
since for the case $\frac{\partial^\top V}{\partial \bm x^\top} \bm G(\bm x) \neq 0$ it is always possible to find an input $\uss$ fulfilling \eqref{eq:bedingung_clf}. 

It is well known that a CLF ensures the existence of an input $\bm u = - \bm k (\bm x)$ such that the closed-loop system is asymptotically stable \cite{Sontag.1989}.
In particular, provided that a suitable CLF is given, it was shown in \cite{Sackmann2000} how for general input-affine nonlinear systems \eqref{input-affine-dynamics}, the following MOC problem can be solved explicitly:
\begin{prop}\cite{Sackmann2000}\label{prop:mor}
Let $V(\bm x)$ be a CLF of \eqref{input-affine-dynamics}. Then an exact solution of the MOC problem
\begin{subequations}\label{eq:mor}
\begin{IEEEeqnarray}{C"l}
		\min_{\uss}  &\frac{1}{2}\int_{\tini}^{\infty} \Qopt + \uss^\top\Ropt\uss \, \dt \\
		\mathrm{s.t.}  &\dot{V}(\x) = \frac{\partial^\top V}{\partial \bm x^\top} \big(\mv{f}(\x) + \gphs(\bm x) \uss \big).
\end{IEEEeqnarray}
\end{subequations}
with $r(\bm x) \succ 0$, $\bm S(\bm x) \succ 0$
is given by
\begin{align}\label{eq:u_mor}
	\uss^\star = -\Ropt^{-1} \gphs^\top(\bm x) \Vx \la(\bm x)
\end{align}
where
\begin{align}
	\la(\x) &:= \frac{\fl + \sqrt{(\fl)^2 + \Ql \cdot \Sl}}{\Sl} \label{eq:lambda_mor},\\
	\fl(\bm x) &:= \frac{\partial^\top V}{\partial \bm x^\top}\bm f(\bm x) \label{eq:def1_simplicity}\\
	\Sl(\bm x) &:= \frac{\partial^\top V}{\partial \bm x^\top} \gphs(\bm x) \mm{S} \gphs^\top(\bm x)\frac{\partial V}{\partial \bm x} \label{ess-lambda}\\
	\Ql(\bm x) &:= \Qopt.\label{eq:def2_simplicity}
	\end{align}
\end{prop}
\begin{pf}
The proof follows the lines of \cite{Sackmann2000} and is listed here for the sake of completeness.
For optimization problem \eqref{eq:mor}, we get the Hamiltonian
\begin{align}
	\Ham(\x,\uss,\la) =& \frac{1}{2} \Qopt + \frac{1}{2}\uss^\top\Ropt\uss  \nonumber \\
	&+ \la \frac{\partial^\top V}{\partial x^\top} \big(\mv{f}(\x) + \gphs(\bm x) \uss \big) \label{Hamiltonian}
\end{align}
with scalar Lagrange multiplier $\Upsilon$. Application of the control equation 
\begin{align}
\frac{\partial 	\Ham(\x,\uss,\la)}{\partial \uss} \stackrel != \nullvec{}
\end{align}
leads to \eqref{eq:u_mor}. From the HJBE for time-invariant systems 
\begin{equation}\label{eq:HJB1}
	0 = \min_{\uss}\, \Ham(\x,\uss,\la)
\end{equation}  
we get 
\begin{align}
	 0 = &\frac{1}{2}\Qopt + \la \frac{\partial^\top V}{\partial \bm x^\top} \mv{f}(\x) \nonumber \\
	 &-\frac{1}{2} \la^2 \frac{\partial^\top V}{\partial \bm x^\top} \gphs(\bm x) \mm{S}^{-1}(\bm x) \gphs^\top(\bm x) \Vx.\label{eq:hjb-mor}
\end{align}
Since \eqref{eq:hjb-mor} is a quadratic function in $\Upsilon$, it has the explicit solution
\begin{subequations}\label{Upsilon-fallunterscheidung}
\begin{numcases}{\la(\bm x)=}
\frac{\fl \pm \sqrt{(\fl)^2 + \Ql \cdot \Sl}}{\Sl}, & $\frac{\partial^\top V}{\partial \bm x^\top} \bm G(\bm x) \neq \nullvec{}$,\label{Upsilon-fall-oben}\\
-\frac{\Ql}{2 \fl}, &  $\frac{\partial^\top V}{\partial \bm x^\top} \bm G(\bm x) = \nullvec{}$\label{Upsilon-fall-unten},
\end{numcases}
\end{subequations}
where \eqref{eq:def1_simplicity}--\eqref{eq:def2_simplicity} are used for compactness of notation.
Note that the ``$+$'' solution in \eqref{Upsilon-fall-oben} implies $\Upsilon(\bm x) > 0$, whereas the
``$-$'' solution in \eqref{Upsilon-fall-oben} is discarded (cf. \eqref{eq:lambda_mor}) since it implies $\Upsilon(\bm x) < 0$ which always leads to an unstable solution.
Moreover, we note that $\frac{\partial^\top V}{\partial \bm x^\top} \bm G(\bm x)=\nullvec{}$ implies $\frac{\partial^\top V}{\partial \bm x^\top} \bm f(\bm x) < 0$, since $V(\bm x)$ is a CLF.
It can be shown \cite[pp.~88,~186]{SackmannDiss} that the Lagrange multiplier $\Upsilon(\bm x)$ in \eqref{Upsilon-fallunterscheidung} is continuous even for the case $\frac{\partial^\top V}{\partial \bm x^\top} \bm G(\bm x) = \nullvec{}$. Hence, it can be fully described with \eqref{eq:lambda_mor} and the distinction of \eqref{Upsilon-fallunterscheidung} is not necessary.
\qed \end{pf}
\begin{rem}
For each equilibrium of the closed-loop system \eqref{input-affine-dynamics}, \eqref{eq:u_mor}--\eqref{eq:def2_simplicity}, asymptotic stability and even hyperstability can be proven \cite[p.~91ff.]{SackmannDiss}. For this purpose, $V(\bm x)$ serves as a Lyapunov function, and with the help of the HJBE \eqref{eq:HJB1}
\begin{align}\label{eq:HJB}
	\min_{\uss}\, \Ham(\x,\uss,\la) &=  \frac{1}{2} r(\x) + \frac{1}{2}(\uss^{\star})^\top\Ropt\uss^\star + \la \dot{V}(\bm x) \nonumber \\
	&= 0,
\end{align}
the negative definiteness of $\dot{V}(\bm x)$ can be ensured, since  
\begin{align}
	 \forall \, \x \neq \nullvec{n}: &&\dot{V}(\bm x) &= -\frac{1}{2\Upsilon(\bm x)}\big( \Qopt + (\uss^{\star})^\top \Ropt \uss^\star \big) \nonumber  \\
	 && & < -\frac{1}{2\Upsilon(\bm x)} \Qopt \leq 0.   \label{eq:stabilitaet_mor}
\end{align}
\end{rem} 
%
\subsection{Control-Lyapunov Functions for Port-Hamiltonian Systems}\label{ch:PHSCLF}
In the following, we consider the case that \eqref{input-affine-dynamics} is a nonlinear or linear ISO-PHS and investigate under which conditions the Hamiltonian $H(\bm x)$ is a CLF. 
\begin{prop}\label{prop:HistCLF}
	Consider an ISO-PHS as in \eqref{eq:systemdynamik_single_player},\eqref{eq:output}. The Hamiltonian \( H(\x) \) is a CLF for \eqref{eq:systemdynamik_single_player} if and only if
	\begin{align}\label{eq:HistCLF2}
		\forall \x \in \Zx: \quad \dHtop \Rphs(\bm x) \dH > 0
	\end{align}
	with $ \Zx := \{ \x \in \Rx \, | \, \gphs^\top(\bm x) \dH = \mv{0},\; \bm x \neq \nullvec{} \} $.
\end{prop}  
\begin{pf}
	Since $\bm J(\bm x)$ is skew-symmetric and $H(\x)$ is a CLF, it holds by definition for all $\bm x \neq \nullvec{}$ (cf. \eqref{eq:bedingung_clf}) that
	\begin{align} \label{eq:HistCLF}
	&	\inf_{\uss} \{\dot{H}(\x) \} \nonumber \\
	& = \inf_{\uss} \bigg\{ \dHtop \big( (\Jphs(\bm x) - \Rphs(\bm x))\dH + \gphs (\bm x)\uss \big) \bigg\} \nonumber \\
		& = \inf_{\uss} \bigg\{ -\dHtop \Rphs(\bm x)\dH + \dHtop\gphs(\bm x)\uss  \bigg\} \nonumber \\
		&< 0.
	\end{align} 	
	If $ \dHtop\gphs(\bm x) \neq \nullvec{m} $, then there always exists
	an input $\bm u'$ such that $ \dHtop\gphs(\bm x) \bm u' < \nullvec{m} $ is fulfilled.
		If $ \dHtop\gphs(\bm x) = \nullvec{m}$, then the first term inside the brackets in \eqref{eq:HistCLF} needs to be negative whenever $\bm x \neq \nullvec{}$, i.e.
	\begin{align}
	 \dHtop \gphs(\bm x) = \nullvec{m} \quad \Rightarrow \quad -\dHtop \Rphs(\bm x) \dH < 0,
	\end{align}
	which is equivalent to \eqref{eq:HistCLF2}. 
	
	Conversely, if \eqref{eq:HistCLF2} is satisfied, the definition of a CLF is automatically fulfilled since $\Zx$ is the manifold where $\frac{\partial^\top V}{\partial \bm x^\top} \gphs(\x) = \nullvec{m}$ in Def. \ref{def:clf} and hence \eqref{eq:HistCLF2} is equivalent to \eqref{def:clf:affine}.
\qed \end{pf}
\begin{cor}
Consider an ISO-PHS system as in \eqref{eq:systemdynamik_single_player},\eqref{eq:output}. The Hamiltonian \( H(\x) \) is a CLF for \eqref{eq:systemdynamik_single_player} if and only if \eqref{eq:systemdynamik_single_player},\eqref{eq:output} is zero-state detectable.
\end{cor}
\begin{pf}
	Since the ISO-PHS \eqref{eq:systemdynamik_single_player} is equipped with the passive output \eqref{eq:output},
	$\mathcal X_G = \{\bm x \in \mathds R^n | \bm y(\bm x)=0\}$ and thus
condition \eqref{eq:HistCLF2} is identical with the definition of a zero-state detectable input-affine nonlinear system \eqref{input-affine-dynamics} given in \cite[p.~47]{vanderSchaft.2017}. 
\end{pf}
\begin{cor}\label{cor:linearPHSCLF}
	Consider an ISO-PHS as in \eqref{eq:systemdynamik_single_player},\eqref{eq:output}. The Hamiltonian \( H(\x) \) is a CLF, if
	\begin{align}
		\forall \x \in \Rx: \quad\mathrm{rank}(\Rphs(\bm x)) = n.
	\end{align}
\end{cor}
\begin{pf}
	Trivially, if $\Rphs(\bm x)$ has full rank, then $ \dHtop \Rphs(\bm x) \dH$ is positive and \eqref{eq:HistCLF2} is fulfilled for all $\x \in \Rx$, which also implies for all $\x \in \Zx$. 
\qed \end{pf}
Next, we introduce a necessary and sufficient condition under which the Hamiltonian of a \emph{linear} ISO-PHS is a CLF:
\begin{prop}\label{cor:HistCLF_linear}
	Consider the linear ISO-PHS dynamics
	\begin{align}
	\dot{\bm x} = \left(\bm J(\bm x) - \bm R(\bm x)\right)\frac{\partial H}{\partial \bm x} + \bm G(\bm x) \bm u
	\end{align}
	with $\bm J = - \bm J^\top$, $\bm R \succeq 0$, $H(\x) = \frac 12 \bm x^\top \bm Q \bm x$ and $\bm Q \succ 0$. Then $H(\bm x)$ is a CLF if and only if
	\begin{align}\label{eq:HistCLF_linear}
		\ker\{  \gphs^\top \mm{Q}\} \, \cap \, \ker\{\mm{Q}^\top \Rphs \mm{Q}\} = \emptyset.
	\end{align}
\end{prop}
\begin{pf}
	For linear ISO-PHSs, the set  $\Zx = \{ \x \in \Rx \, | \, \gphs^\top \dH = \mv{0} \} $ is equivalent to the kernel of $\gphs^\top \mm{Q}$ since
	\begin{align}
	\forall \; \x \in \ker\{ \gphs^\top \mm{Q}\}: && 	\gphs^\top \dH =  \gphs^\top \mm{Q} \x = \nullvec{n}.
	\end{align}
	Thus following Proposition~\ref{prop:HistCLF}, 
	\begin{align}\label{eq:bedingunganRlinear}
	\forall \x \in  \ker\{ \gphs^\top \mm{Q}\}: && 	-\dHtop \Rphs \dH = -\mm{Q}\x^\top \Rphs \mm{Q} \x < 0 
	\end{align}
	has to be satisfied. Firstly, it is important to note that \eqref{eq:bedingunganRlinear} is always $\leq 0$ for $\x \in \Rx$, since the product $\mm{Q}^\top \Rphs \mm{Q}$ of positive-semidefinite matrices is again positive-semidefinite \cite[p.~431]{horn2012matrix}. Secondly, with Lemma~\ref{lemma:eigenwert_und_kern} (see Appendix \ref{app:Kernlemma}), follows that the equality $\x^\top \mm{Q}^\top \Rphs \mm{Q} \x = 0$ holds if and only if $\x \in \ker\{\mm{Q}^\top \Rphs \mm{Q}\}$.
	Consequently, \eqref{eq:bedingunganRlinear} holds if and only if
	 $\ker\{  \gphs^\top \mm{Q}\}$ and $\ker\{\mm{Q}^\top \Rphs \mm{Q}\}$ are disjoint. 
\qed \end{pf}
\subsection{Relationship between Modified Optimal Control and Optimal Control}\label{ch:MOROC}
Taking Proposition~\ref{prop:HistCLF} into account, we can apply MOC to ISO-PHSs in a straightforward manner. However, the question arises to what extent the arising controller \eqref{eq:u_mor}--\eqref{eq:def2_simplicity} is optimal with respect to the original optimal control problem \eqref{original-OP}. 
\begin{prop}\label{prop:lambda1}
	The (modified optimal) controller \eqref{eq:u_mor}--\eqref{eq:def2_simplicity} is optimal with respect to \eqref{original-OP} if
	\begin{align}
		\forall \x \in \mathds R^n: &&  \la(\x) = 1. \label{UpsilonGleichEins}
	\end{align}
\end{prop}
\begin{pf}
If $\Sl \neq 0$ then it follows
 from \eqref{ess-lambda} that $\bm G^\top \frac{\partial V}{\partial \bm x} \neq 0$. Thus \eqref{UpsilonGleichEins} is equivalent to
 	\begin{align}\label{eq:lambda_kurz}
		\la(\x)= \frac{\fl + \sqrt{(\fl)^2 + \Ql \cdot \Sl}}{\Sl} \stackrel{!}{=} 1.
	\end{align}
	From \eqref{eq:lambda_kurz} we obtain	
	\begin{align}
	\Sl - 2\fl = \Ql. \label{miniQ}
	\end{align}
	Substitution of \eqref{eq:def1_simplicity}-\eqref{eq:def2_simplicity} into \eqref{miniQ} leads to
	\begin{align}
		\frac{1}{2} \Qopt - \frac{1}{2} \Sx + \vtf = 0, \label{eq:proof-lambda1-umformung4}
	\end{align}
	which is exactly the HJBE for time-invariant systems \eqref{eq:HJB}. The function $V(\x)$ solving \eqref{eq:proof-lambda1-umformung4} is the value function $V^\star(\x)$. This means, if we achieve to find a CLF for which $\la(\x) = 1$ holds for all $\bm x \in \mathds R^n$, this CLF is also the value function $V^\star$ and thus \eqref{eq:u_mor} is an optimal control input for \eqref{original-OP}. 
	
	If $\Sl = 0$, then $\bm G^\top(\bm x) \frac{\partial V}{\partial \bm x} = 0$. Thus \eqref{UpsilonGleichEins} is equivalent to
 	\begin{align}\label{eq:lambda-ganzkurz}
\la(\x)= - \frac{\Ql}{2 \fl} \stackrel{!}{=} 1,
\end{align}	
	which also leads to \eqref{miniQ}.
\qed \end{pf}
A less restrictive requirement can be derived by allowing $\Upsilon (\bm x)$ to be an arbitrary, but fixed positive value:
\begin{cor}\label{valuefunctionsuper}
	The (modified optimal) controller \eqref{eq:u_mor}--\eqref{eq:def2_simplicity} is optimal with respect to \eqref{original-OP} if
	\begin{align}
	\forall \x \in \mathds R^n:  && \la(\x) = c, && c \in \R_{>0}.
	\end{align}
	In this case, the value function $V^\star(\bm x)$ is a $c$-multiple of the CLF $V(\bm x)$, i.e. $V^\star(\bm x) = c \cdot V(\bm x)$.
\end{cor}
\begin{pf} Following the same procedure as in Proposition~\ref{prop:lambda1}, the HJBE \eqref{eq:proof-lambda1-umformung4} becomes
	\begin{align}
		\frac{1}{2} \Qopt - \frac{1}{2} c^2 \Sx + c\vtf = 0.
	\end{align}	
	Accordingly, it follows that $c \cdot V(\x)$ is the value function. 
\qed \end{pf}
\begin{rem}\label{rem:suboptimal}
From Corollary \ref{valuefunctionsuper} we can conclude that unless $\Upsilon(\bm x)$ converges to a constant value and remains constant even after a disturbance, the chosen CLF cannot be equivalent to the value function. Consequently the resulting controller \eqref{eq:u_mor}--\eqref{eq:def2_simplicity} is not optimal with respect to \eqref{original-OP}. The fluctuation of $\Upsilon(t)$ over time can therefore be interpreted as an \emph{indicator of suboptimality} with respect to \eqref{original-OP}.
\end{rem}
\section{From Modified to Optimal Control}
Since $H(\bm x)$ is in general not equal to the value function and hence the condition of Proposition \ref{prop:lambda1} is not fulfilled by using $V(\bm x)=H(\bm x)$, the modified optimal controller presented in the previous chapter is not optimal with respect to the original problem \eqref{original-OP}.
For this purpose, in Section 3.1, we present an extended CLF $V(\bm x, \bm w)$ as a
linear combination of $H(\bm x)$ and a weighted set of basis functions.
In Section 3.2, we propose a gradient-based adaptation strategy of the weighting factors in $V(\bm x, \bm w)$ such that $V(\bm x, \bm w)$ fulfills the condition of Proposition \ref{prop:lambda1}.
In Section 3.3, we deploy an adaptive optimal controller for \eqref{original-OP} based on $V(\bm x, \bm w)$ and MOC and show that the equilibrium of the closed-loop system is (asymptotically) stable.
\subsection{Extended Control-Lyapunov Function}
The extended CLF is composed of $H(\bm x)$ and a weighted sum of $r \in \mathds N$ basis functions:
\begin{align}
V(\bm x,\bm w)=H(\bm x) + \bm w^\top \bm \Phi(\bm x), \label{ansatz}
\end{align}
where $\bm \Phi: \mathds R^n \to \mathds R^r$ and $\bm w \in \mathds R^r$. 
We assume that the basis functions $\bm \Phi(\bm x)$ are $C^2$ and ``properly chosen'' in the sense that
the actual value function $V^\star$ can be parameterized via $\bm \Phi(\bm x)$ and an optimal weighting vector $\bm w^\star$:
\begin{assum}\label{ass:optimal}
	\begin{align}
	\exists \bm w^\star \in \mathds R^r: && V^\star (\bm x) = (\bm w^\star)^\top \bm \Phi(\bm x)
	\end{align}
\end{assum}
This assumption is admissible, if the number of basis functions is large (see e.g. \cite[p.~881]{Vamvoudakis.2010},\cite[p.~1020f.]{Liu.2014}), and later allows to characterize the deviation from the optimal solution by the distance between $\bm w(t)$ and $\bm w^\star$.

With \eqref{ansatz}, we obtain
\begin{align}
	\frac{\partial V(\bm x, \bm w)}{\partial \x} = \dH + \frac{\partial^\top \bm \Phi}{\partial \bm x^\top}\w
	\end{align}
and accordingly the MOC law \eqref{eq:u_mor} reads
\begin{align}
\uss^\star = -\Ropt^{-1} \gphs^\top \Vx \la(\bm x, \bm w), \label{mor-mit-adaption}
\end{align}
where
	\begin{align}
	\la(\x, \bm w) &= \frac{\fl' + \sqrt{(\fl')^2 + \Ql' \cdot \Sl'}}{\Sl'}, \label{neueKurzsymboleAnfang}\\
		\fl'(\bm x)& = \frac{\partial ^\top H}{\partial \bm x^\top}(\bm J(\bm x)-\bm R(\bm x))\dH \nonumber \\
		& + \w^\top \frac{\partial \bm \Phi}{\partial \bm x} (\bm J(\bm x)-\bm R(\bm x))\dH, \\
		\Sl'(\bm x) &= \w^\top \frac{\partial \bm \Phi}{\partial \bm x} \mm{K}(\bm x) \frac{\partial^\top \bm \Phi}{\partial \bm x^\top} \w + 2 \, \frac{\partial ^\top H}{\partial \bm x^\top}\mm{K}(\bm x) \frac{\partial^\top \bm \Phi}{\partial \bm x^\top} \w \nonumber \\  \label{Sl-strich}
		& + \frac{\partial ^\top H}{\partial \bm x^\top} \mm{K}(\bm x)\dH, \\
		\Ql'(\bm x)&= \Qopt, \label{ZusammenhangQundr}\\
		\bm K(\bm x) &= \bm G(\bm x) \bm S(\bm x) \bm G^\top(\bm x). \label{neueKurzsymboleEnde}
	\end{align}
Employing the same reasoning as in Proposition~\ref{prop:lambda1}, we will study in more detail how to check whether a given CLF \eqref{ansatz} is equivalent to the value function $V^\star$.
\begin{prop}\label{cond:quadrik}
	Let $V(\bm x) = H(\bm x) + (\bm w^\diamond)^\top \bm \Phi(\bm x)$ be a given CLF with $\bm w^\diamond \in \mathds R^r$. Then
	$V(\bm x)$ is equivalent to the value function $V^\star(\bm x)$ of \eqref{original-OP} if and only if 
	\begin{align}
	\forall\, \x \in \Rx: && (\bm x, \w^\diamond) \in \Qu(\x, \bm w), \label{bedingung:prop:lambda1}
	\end{align}
	where
	\begin{align}
	\Qu(\x. \bm w) &= \left\{	 (\bm x, \bm w) \in \mathds R^n \times \mathds R^r:\right. \nonumber \\
	& \quad \quad \left.\w^\top \bm A(\bm x) \w 	+ \bm a^\top(\bm x)\w    +  a(\bm x) = 0 \right\} \label{quadrik}	\\	
		 \bm A(\bm x)&= \frac{\partial \bm \Phi}{\partial \bm x} \bm K(\bm x) \frac{\partial^\top \bm \Phi}{\partial \bm x^\top} \label{grosses-A} \\
	\bm a(\bm x)&= 2\bigg(\frac{\partial ^\top H}{\partial \bm x^\top}\bm K(\bm x) \frac{\partial^\top \bm \Phi}{\partial \bm x^\top}  \nonumber \\
	 &- \frac{\partial ^\top H}{\partial \bm x^\top}(\bm J(\bm x)-\bm R(\bm x))\frac{\partial^\top \bm \Phi}{\partial \bm x^\top} \bigg)  \\ 
	a(\bm x)&=\frac{\partial ^\top H}{\partial \bm x^\top} \bm K(\bm x)\dH \nonumber \\
	&- 2 \frac{\partial ^\top H}{\partial \bm x^\top}(\bm J(\bm x)-\bm R(\bm x))\dH - \Qopt \label{grosses-K}
	\end{align}
\end{prop}
\begin{pf}
Let $\bm w^\ast \in \mathds R^r$ be the optimal weighting vector. According to Proposition~\ref{prop:lambda1}, this implies that 
\begin{align}\label{eq:proofNC1}
\forall\, \x \in \Rx: && \la(\x, \bm w^\ast) = 1.
\end{align}
As shown in the proof of Proposition~\ref{prop:lambda1}, condition \eqref{eq:proofNC1} is equivalent to
\begin{align}
\Sl' - 2\fl' - \Ql'
\stackrel != 0. \label{bedingung:neuerBeweis}
\end{align}
Inserting of \eqref{neueKurzsymboleAnfang}--\eqref{neueKurzsymboleEnde} in \eqref{bedingung:neuerBeweis} yields
\begin{align}
& (\w^\star)^\top \frac{\partial \bm \Phi}{\partial^\top \bm x} \mm{K}(\bm x) \frac{\partial \bm \Phi^\top}{\partial \bm x} \w^\star  + 2\bigg(\frac{\partial ^\top H}{\partial \bm x^\top}\mm{K}(\bm x) \frac{\partial^\top \bm \Phi}{\partial \bm x^\top}  \nonumber \\
&- \frac{\partial ^\top H}{\partial \bm x^\top}(\bm J(\bm x)-\bm R(\bm x))\frac{\partial^\top \bm \Phi}{\partial \bm x^\top} \bigg)\w^\star   \nonumber  \\ 
&+  \frac{\partial ^\top H}{\partial \bm x^\top} \mm{K}(\bm x)\dH - 2 \frac{\partial ^\top H}{\partial \bm x^\top}(\bm J(\bm x)-\bm R(\bm x))\dH - \Qopt \nonumber \\
&= 0. \label{eq:quadrik-condition}
\end{align}
With $\mathcal Q(\bm x, \bm w)$ as in \eqref{quadrik}, condition \eqref{eq:quadrik-condition} can be written as 
\begin{align}
\forall \, \x \in \Rx: && (\bm x, \bm w^\star) \in \mathcal Q(\bm x, \bm w). 
\end{align}
Since the given CLF $V(\bm x)$ is equal to $V^\star(\bm x)$ if and only if $\bm w^\diamond = \bm w^\star$, this is equivalent to \eqref{bedingung:prop:lambda1}. \qed
\end{pf}
Since $\mathcal Q(\bm x, \bm w)$ can be written as 
\begin{align}
	\Qu(\x. \bm w) &= \left\{	 (\bm x, \bm w) \in \mathds R^n \times \mathds R^r: Q(\bm x, \bm w) = 0 \right\}	
\end{align}
with
\begin{align}
Q(\x,\w) = \bm w^\top \bm A(\bm x) \bm w + \bm a^\top(\bm x) \bm w + a(\bm x)
\end{align}
being a quadratic function, $\mathcal Q(\bm x, \bm w)$ is a quadric for each fixed $\bm x \in \mathds R^n$. It has two important properties:
Firstly, the shape of the quadric is dependent on $\bm x$. Secondly, according to Proposition~\ref{cond:quadrik}, the optimal weighting vector $\bm w^\star$ is contained in each quadric and thus
\begin{align}
\forall\, \bm x \in \mathds R^n: && Q(\bm x, \bm w^\star) = 0. \label{wast-nullt-Q}
\end{align}
We will exploit both of these facts in Section \ref{thbehnjrez} to derive a gradient descent procedure ensuring convergence to $\bm w^\star$.
%
%
%
%
%
\subsection{Adaptation of the extended Control-Lyapunov Function}\label{thbehnjrez}
As shown in \eqref{wast-nullt-Q}, the optimal weighting vector $\bm w^\star$ is a root of $Q(\bm x, \bm w)$, independent of $\bm x$.
Thus for each arbitrary but fixed $\bm x \in \mathds R^n$, we can characterize $\bm w^\ast$ as the minimizer of an objective function $J_0(\bm x, \bm w)$ with
\begin{align}
J_0(\bm x, \bm w):=\left(Q(\bm x, \bm w)\right)^2. \label{jott-null}
\end{align} 
Moreover, with $\mathcal Q(\bm x, \bm w) = \argmin_{\bm w} \left\{ J_0 (\bm x,\bm w)\right\}$ and due to the fact that $\bm w^\ast$ is contained in each quadric $\mathcal Q(\bm x, \bm w)$, it follows that   
\begin{align}
\{ \bm w^\ast\} \subseteq \bigcap_{\bm x \in \mathds R^n} \argmin_{\bm w} \left\{ J_0 (\bm x,\bm w)\right\}.
\end{align}
However, we note that $J_0$ is in general \emph{not} strictly convex around $\bm w^\star$, which hampers convergence to $\bm w^\star$ and necessitates additional conditions for a sufficient exploration of the state space. To circumvent these, often very hard-to-evaluate requirements, we formulate an extended objective function $J_w(\bm x, \bm w)$ 
providing strict convexity with respect to $\bm w$ in a neighborhood of $\bm w^\star$.
This is outlined in the next proposition.
\begin{prop}\label{prop:konvexitat}
	Let 
	\begin{align}
	Q(\x,\w) = \bm w^\top \bm A(\bm x) \bm w + \bm a^\top(\bm x) \bm w + a(\bm x)
	\end{align}
	with $\bm A(\bm x)$, $\bm a(\bm x)$, $a(\bm x)$ as in \eqref{grosses-A}--\eqref{grosses-K}
	be the corresponding quadratic function to the quadric $\mathcal Q(\bm x, \bm w)$ in \eqref{quadrik}. Then with $\bm c_1, \ldots, \bm c_r \in \mathds R^n$, the extended objective function
	\begin{align}
	J_w(\bm x, \bm w) &= J_0(\bm x+\bm c_1, \bm w) + \cdots + J_0(\bm x+\bm c_r, \bm w) \nonumber  \\
	&=(Q(\x+\mv{c}_1,\w))^2  + \dots +  (Q(\x+\mv{c}_{\AnzahlGewichte},\w))^2 \label{eq:linearkombinationNeu}
	\end{align}
	composed by a linear combination of shifted objective functions $J_0(\bm x, \bm w)$ is (locally) strictly convex in an open neighborhood $\mathcal{N}$ of the optimal weights $\w^\star$ for all $\x \in \Rx$, if and only if the vectors
	\begin{align}\label{eq:vektoren_v}
		\mv{v}_i = 2\mm{A}(\x+ \mv{c}_i) \w^\star + \mv{a}(\x+\mv{c}_i)
	\end{align}
	with $i = 1, \dots, \AnzahlGewichte$ are linearly independent. 
\end{prop}  
\begin{pf}
	Without loss of generality, we choose $\bm c_1 = \nullvec{}$. As strict convexity with respect to $\bm w$ needs to be shown, the Hessians of $J_0$ (see \eqref{jott-null}) and $J_w$ (see{\eqref{eq:linearkombinationNeu}}) are studied. The Hessian of $J_0$ can be written as  
	\begin{align}
		\frac{\partial^2 J_0(\bm x, \bm w)}{\partial \bm w^2}  &=  2\big(2\mm{A}(\x) \w + \mv{a}(\x)\big)\big(2\mm{A}(\x) \w + \mv{a}(\x)\big)^\top \nonumber \\  & +  4\mm{A}(\x)\big(\w^\top \mm{A}(\x) \w + \mv{a}(\x) \w + a(\x)\big).\label{eq:hessian_quadrik}
	\end{align}
	Since the optimal weights $\w^\star$ need to be part of each quadric $\Qu(\x, \bm w)$ regardless of the state $\x$ (see Proposition~\ref{cond:quadrik}), 
	the second summand in \eqref{eq:hessian_quadrik} is equal to zero for $\bm w = \bm w^\star$ and accordingly
	\begin{align}\label{eq:hessian_rank1}
		\frac{\partial^2 J_0(\x,\w)}{\partial \w^2}\bigg\vert_{\w = \w^\star} = 2 \mv{v}_1 {\mv{v}_1^\top}, 
	\end{align} 		
	with $\mv{v}_1 = 2\mm{A}(\x) \w^\star + \mv{a}(\x)$.
	Since the Hessian \eqref{eq:hessian_rank1} is only composed by the multiplication of two vectors, which yields a matrix with identical but scaled row vectors, it is positive-semidefinite and has rank one.
	
	For the Hessian of the shifted objective function $J_0(\bm x + \bm c_2 , \bm w)$, we obtain in a similar manner
	\begin{align} \label{eq:hessian_rank1-shift}
		\frac{\partial^2 J_0(\x + \mv{c}_2, \w)}{\partial \w^2}\bigg\vert_{\w = \w^\star} = 2 \mv{v}_2 {\mv{v}_2^\top}
	\end{align}
	with $\mv{v}_2 = 2\mm{A}(\x + \mv{c}_2) \w^\star + \mv{a}(\x + \mv{c}_2)$. Note that the rank of the matrix \eqref{eq:hessian_rank1-shift} is one regardless of the shifting $\mv{c}_2$. 
	
	The linear combination $J_{1}(\bm x, \bm w) := J_0(\x,\w) + J_0(\x + \mv{c}_2,\w)$ leads to the Hessian
	\begin{align}
		\frac{\partial^2 J_1}{\partial \w^2}\bigg\vert_{\w = \w^\star} = 2 \mv{v}_1 {\mv{v}_1^\top} + 2 \mv{v}_2 {\mv{v}_2^\top}. \label{zweierkombination}
	\end{align}	
	The same applies to linear combinations with more than two summands due to the linearity property of differentiation. 
	
	We can see that \eqref{zweierkombination} has a maximum rank of two. As strict convexity of $\bm J_w$ is required, full rank $r$ needs to be satisfied for the Hessian of $J_w(\bm x, \bm w)$ at $\bm w=\bm w^\star$. Thus the question arises, in which case the increase of summands implies a rank increase of the Hessian. 
	Each matrix in \eqref{eq:hessian_rank1} or \eqref{eq:hessian_rank1-shift} describes a linear map $\Rw \rightarrow \Rw$ of rank one and its image is a subspace of $\Rw$ of dimension one. If the image of $\frac{\partial^2 J_1}{\partial \w^2}$ in \eqref{zweierkombination} is of dimension two, the rank automatically increases, since 
	\begin{align}
		\dim(\text{im}(\mm{M})) = \rank(\mm{M}).
	\end{align} 
	for an arbitrary matrix $\mm{M}$. With the dimension formula for the sum of subspaces \cite[p.~47]{Bosch.2014}, it follows for two arbitrary subspaces $\mathcal U_1$ and $\mathcal U_2$
	\begin{align}
		\dim(\mathcal U_1 + \mathcal U_2) = \dim(\mathcal U_1) + \dim(\mathcal U_2) - \dim(\mathcal U_1 \cap \mathcal U_2).
	\end{align}
	By setting $\mathcal U_1 = \text{im}(\mv{v}_1 {\mv{v}_1}^\top)$ and $\mathcal U_2 = \text{im}(\mv{v}_2 {\mv{v}_2}^\top)$ it follows that only if $\dim(\mathcal U_1 \cap \mathcal U_2) = 0$, the sum of the two matrices $\mv{v}_1 {\mv{v}_1}^\top$ and $\mv{v}_2 {\mv{v}_2}^\top$ leads to an increase of rank. The rank-one matrix $\mv{v}_l {\mv{v}_l}^\top$ is formed by weighted rows of ${\mv{v}_l}^\top$ with the respective components $v_{l_i}$, $i = 1,\dots, r$ of ${\mv{v}_l}$ 
	\begin{align}
		\mv{v}_l {\mv{v}_l}^\top = \begin{bmatrix}v_{l_1} \cdot {\mv{v}_l}^\top \\ \vdots \\ v_{l_r} \cdot {\mv{v}_l}^\top \end{bmatrix}, 
	\end{align}
	and hence its image spans the subspace
	\begin{align}
		\mathcal U_l = \text{im}(\mv{v}_l {\mv{v}_l}^\top) = \{\mu \cdot \mv{v}_l \mid \mu \in \R \}.
	\end{align}
	With regard to $J_1$ in \eqref{zweierkombination} we see that in order to let both sets $\mathcal U_1$ and $\mathcal U_2$ be disjunct, the linear independence of both vectors $\bm v_1$ and $\bm v_2$ is necessary. Graphically, the subspace $\text{im}(\mv{v} {\mv{v}}^\top)$ is a straight line in $\Rw$, and linear independence leads to non coinciding straight lines such that $\dim(\mathcal U_1 \cap \mathcal U_2) = 0$. 
	
	Applied to $J_w$ it becomes clear that \eqref{eq:linearkombinationNeu} has a Hessian with full rank if the vectors $\mv{v}$ in \eqref{eq:vektoren_v} with $l= 1, \ldots, r$ are linearly independent. Hence, each vector $\mv{v}_l$ induces a matrix implying an increase of one for the rank of the Hessian $\frac{\partial^2 J_w}{\partial \w^2}$ which leads to full rank and thus positive definiteness of $\frac{\partial^2 J_w}{\partial \w^2}$. Moreover, since $J_w$ is a $C^2$ function, positive definiteness of the Hessian is preserved for all $\bm w \in \mathcal N$ around $\bm w^\ast$.
\qed \end{pf}
\begin{rem}
	Proposition \ref{prop:konvexitat} states that the minimum number of summands in $J_w$ (see \eqref{eq:linearkombinationNeu}) is $r$ to achieve strict convexity. From a practical point of view, however, it is feasible and preferable to compose $J_w$ from more than $r$ summands in order to enhance convergence properties.
	\end{rem}
Due to the fact that \eqref{wast-nullt-Q} is fulfilled for all $\bm x \in \mathds R^n$, the optimal weighing factor $\bm w^\star$  can be characterized by the strictly convex optimization problem
\begin{align}
\bm w^\star = \argmin_{\bm w} \left\{ J_w (\bm x,\bm w)\right\}.\label{eq:linearkombination}
\end{align}
Thus, a given weighting factor $\bm w$ can be adapted by the gradient descent procedure
\begin{align}
\dot{\bm w} = - \alpha \cdot \frac{\partial J_w(\bm x, \bm w)}{\partial \bm w}
 \label{adaptationNeu}
\end{align}
with learning rate $\alpha >0$.
\subsection{Stability of the Closed-Loop System}
\begin{figure}
	\centering
	\includegraphics[height=0.42\textwidth]{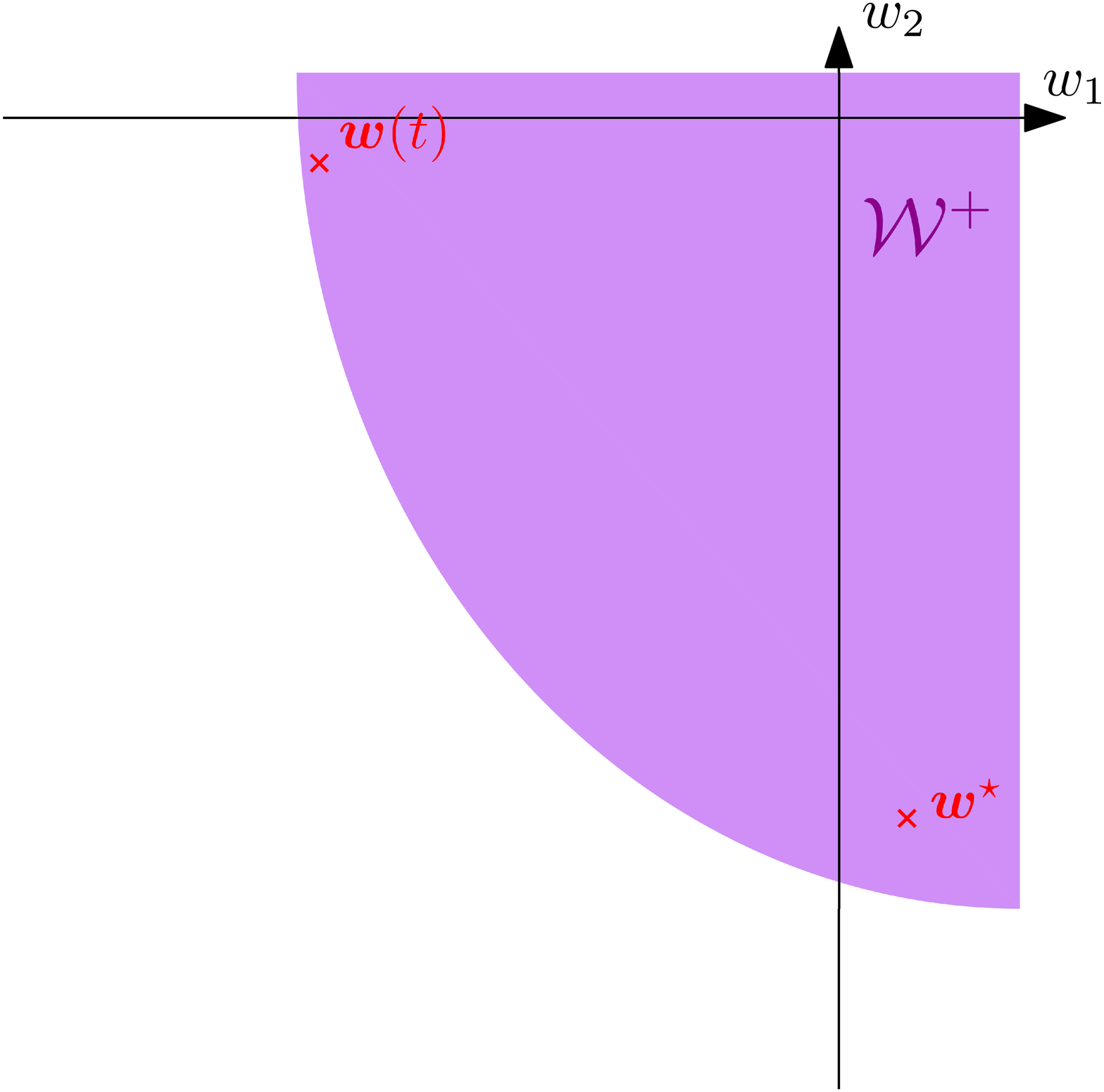}
	\caption{Illustration of the set $\mathcal W^+$ with $\bm w^\star \in \mathrm{int}\, \mathcal W^+$ and $\nullvec{} \in \mathrm{int}\, \mathcal W^+$.}
	\label{fig:Illustration}
\end{figure}

\begin{figure}
	\centering
	\includegraphics[height=0.42\textwidth]{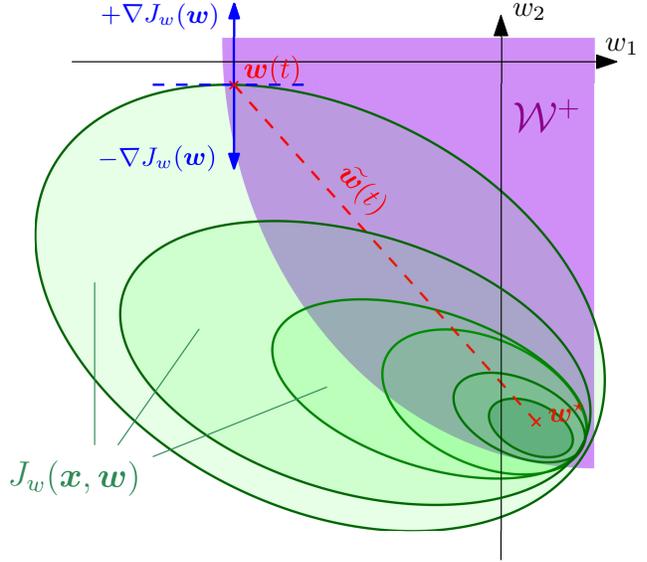}
	\caption{Contour plot of $J_w(\bm x, \bm w)$ for fixed $\bm x$.}
	\label{fig:toWAst}
\end{figure}

We are now ready to formulate an explicit control law which solves the original optimal control problem \eqref{original-OP}.
With the open-loop ISO-PHS \eqref{eq:systemdynamik_single_player}, the extended CLF \eqref{ansatz}, the MOC law \eqref{mor-mit-adaption}, the adaptation procedure \eqref{adaptationNeu}, and the shorthand notation \eqref{neueKurzsymboleAnfang}--\eqref{neueKurzsymboleEnde}, we get the following closed-loop system:
\begin{subequations}\label{closed-loop}
\begin{align}
	\displaystyle \dot{\bm x} &= \displaystyle (\bm J(\bm x) - \bm R(\bm x))\frac{\partial H(\bm x)}{\partial \bm x} + \bm G(\bm x) \bm u^\star \label{closed-loop-strecke} \\
	\displaystyle \bm u^\star & = \displaystyle \bm S^{-1}(\bm x) \bm G^\top(\bm x) \frac{\partial V}{\partial \bm x} \cdot \frac{\fl' + \sqrt{(\fl')^2 + \Ql' \cdot \Sl'}}{\Sl'}  \label{eq:closed_loop_eingang}\\
	\displaystyle V(\bm x, \bm w) &= \bm H(\bm x) + \bm w^\top \bm \Phi(\bm x) \label{vau} \\
	\displaystyle \dot{\bm w} &= \displaystyle- \alpha \cdot \frac{\partial J_w(\bm x, \bm w)}{\partial \bm w} \label{adaptation}\\
	\bm x_0 &= \bm x_0 \\
	\bm w_0 &= \nullvec{r}
\end{align}
\end{subequations}
To perform a stability analysis of the equilibrium $(\nullvec{},\bm w^\star)$ of \eqref{closed-loop}, we use $V(\bm x, \bm w)$ as a Lyapunov function candidate and prove that 
\begin{align}
\forall \; (\bm x, \bm w) \neq (\nullvec{}, \bm w^\star): && V(\bm x, \bm w) > 0 \label{rr} \\
\forall \; (\bm x, \bm w) \in \mathds R^n \times \mathds R^r: && \dot V(\bm x, \bm w) \leq  0 \label{ss}
\end{align} 
While the proof of \eqref{ss} is straightforward (see Proposition~\ref{dot-V-decreasing}), statement \eqref{rr} (see Proposition~\ref{prop:Vposdefinit}) requires some additional preparatory work.
\begin{prop}\label{dot-V-decreasing}
	Consider the closed-loop system \eqref{closed-loop} starting at $(\bm x_0, \bm w_0) \in \mathds R^n \times \mathds R^r$. Then
	\begin{align}
	\forall \; (\bm x, \bm w) \in \mathds R^n \times \mathds R^r: && \dot V(\bm x, \bm w) \leq  0,
	\end{align}
	i.e. $V(\bm x, \bm w)$ decreases monotonically over time.
\end{prop}
\begin{pf}
	Applying the chain rule to \eqref{vau} and inserting \eqref{closed-loop-strecke}, \eqref{eq:closed_loop_eingang}, \eqref{adaptation}, we get
	\begin{align}
	\dot V&(\bm x, \bm w) \nonumber \\
	=& \frac{\partial ^\top V(\bm x, \bm w)}{\partial ^\top \bm x} \dot{\bm x} + \frac{\partial ^\top V(\bm x, \bm w)}{\partial ^\top \bm w} \dot{\bm w} \nonumber \\
	=& \frac{\partial ^\top V(\bm x, \bm w)}{\partial ^\top \bm x} (\bm J(\bm x) - \bm R(\bm x)) \frac{\partial H}{\partial x} \nonumber \\
	&- \frac{\partial ^\top V(\bm x, \bm w)}{\partial ^\top \bm x} \bm G(\bm x) \bm S^{-1}(\bm x) \bm G^\top(\bm x) \frac{\partial V(\bm x, \bm w)}{\partial\bm x} \bm u \nonumber \\
	&- \frac{\partial ^\top V(\bm x, \bm w)}{\partial ^\top \bm w} \alpha  \frac{\partial  J_w(\bm x, \bm w)}{\partial \bm w} \\
	=& \frac{\partial ^\top V(\bm x, \bm w)}{\partial ^\top \bm x} (\bm J(\bm x) - \bm R(\bm x)) \frac{\partial H}{\partial x}  \nonumber \\
		&- \frac{\partial ^\top V(\bm x, \bm w)}{\partial ^\top \bm w} \alpha  \frac{\partial  J_w(\bm x, \bm w)}{\partial \bm w} \nonumber \\
	&	- \Sl' \cdot \frac{\fl' + \sqrt{(\fl')^2 + \Ql' \cdot \Sl'}}{\Sl'} \nonumber \\
		=& \fl' - \fl' -  \sqrt{(\fl')^2 + \Ql' \cdot \Sl'} \nonumber \\
		=& -  \sqrt{(\fl')^2 + \Ql' \cdot \Sl'}. \label{dot_v_result}
	\end{align}
		With \eqref{Sl-strich} in \eqref{dot_v_result}, we get
		\begin{align}
			\Sl' =& \frac{\partial^\top V(\bm x, \bm w)}{\partial \bm x^\top} \gphs(\bm x) \mm{S} \gphs^\top(\bm x)\frac{\partial V(\bm x, \bm w)}{\partial \bm x} \label{FormelFuerSUpsilon}
		\end{align}
		and due to the fact that $\bm S(\bm x) \succ 0$, it follows that $\Sl' \geq 0$. Moreover, $\Ql'\geq 0$ holds per definition (see \eqref{ZusammenhangQundr}). Hence $\dot V(\bm x, \bm w)$ is nonpositive for all $(\bm x, \bm w) \in \mathds R^n \times \mathds R^r$.
\qed \end{pf}
Now the positive definiteness of $V=H(\bm x) + \bm w^\top \bm \Phi(\bm x)$ has to be evaluated. 
Despite the fact that $H(\bm x)$ is positive-definite per definition, $V(\bm x, \bm w)$ may be nonpositive if ``$+\bm w^\top \bm \Phi(\bm x)$'' is negative for some $(\bm x, \bm w) \in \mathds R^n \times \mathds R^r$. In particular, we thus have to prove that $V(\bm x(t), \bm w(t))$ is still positive for all $t \geq t_0$ where $\bm x \neq \nullvec{}$.
For a closer look at this question, let $\mathcal W^+ \subseteq \mathds R^r$ denote the set of $\bm w \in \mathds R^r$ where $V(\bm x, \bm w)>0$ is fulfilled for \emph{all} $\bm x \neq \nullvec{}$:
\begin{align}
\mathcal W^+ = \left\{\bm w \in \mathds R^r | V(\bm x, \bm w) > 0 \; \forall \bm x \neq  \nullvec{} \right\}
\end{align}
We first have to prove several properties of $\mathcal W^+$ to conclude that each trajectory $\bm x \in \mathds R^n$ of \eqref{closed-loop} starting at $(\bm x_0, \nullvec{})$ will converge to the set $\mathcal W^+$.
\begin{lem}\label{wastDrinne}
	$\bm w^\star \in \mathcal W^+$ and $\nullvec{} \in \mathcal W^+$.
\end{lem}
\begin{pf}
	Since $V(\bm x, \bm w^\star)=V^\star (\bm x) = H(\bm x) + (\bm w^\star)^\top \bm \Phi(\bm x)$ is the value function according to Assumption \ref{ass:optimal}, $V^\star(\bm x)\succ 0$ holds per definition and thus $\bm w^\star \in \mathcal W^+$. 
	
	For $\bm w = \nullvec{}$, it trivially holds that $V(\bm x)=H(\bm x) \succ 0$, since positive definiteness of $H(\bm x)$ is fulfilled per definition. Consequently, $\nullvec{} \in \mathcal W^+$.
\qed \end{pf}
As $\mathcal W^+$ is an open set (see Lemma \ref{lemma:wPlusOpen} in Appendix \ref{ch:openLemma}), we conclude that $\nullvec{} \in \mathrm{int}\, \mathcal W^+$ and $\bm w^\star \in \mathrm{int}\, \mathcal W^+$, which is illustrated in Fig. \ref{fig:Illustration}.
Note that this fact does not imply that $\bm w(t) \in \mathcal W^+$ holds for \emph{all} $t$, see Fig. \ref{fig:toWAst} for an illustrative example: We can see the contour plot of $J_w(\bm x, \bm w)$ for a fixed $\bm x$. Of course, $\bm w^\star = \argmin J_w(\bm x, \bm w)$. 
However, depending on the shape of $J_w$, it may be possible that the descent direction $-\nabla J_w$ is pointing out of $\mathcal W^+$, which yields $\bm w(t') \notin \mathcal W^+$ for some $t'>t$.

Despite the fact that $\bm w(t)$ may be temporarily outside of $\mathcal W^+$, we will prove now that for a sufficiently large but finite $T\geq t_0$, $\bm w(t)$ always lies within $\mathcal W^+$. 
This is stated in Proposition~\ref{prop:irgendwann-in-W-plus} by making use of Proposition~\ref{streng-monoton-fallend}.
\begin{prop}\label{streng-monoton-fallend}
	Let $V^\star(\bm x) = H(\bm x) + (\bm w^{\star})^\top \bm \Phi(\bm x)$ be the value function of optimization problem \eqref{original-OP} and let the conditions of Proposition~\ref{prop:konvexitat} hold with $\bm w_0 \in \mathcal N$. Then, the trajectory of $\bm w(t)$ fulfills
	\begin{align}
	\lim_{t \to \infty} \left \| \bm w(t) - \bm w^\star \right\|_2 = 0, \label{w-trajectory}
	\end{align}  
	i.e. all weighting factors $\bm w(t)$ asymptotically converge to the optimal ones.
\end{prop}
\begin{pf}
	If the conditions of Proposition~\ref{prop:konvexitat} hold with $\bm w_0 \in \mathcal N$, then $J_w(\bm x,\bm w)$ is strictly convex with respect to $\bm w$ in an open neighborhood $\mathcal N$ of the optimizer $\bm w^\star$ for each arbitrary but fixed $\bm x$, i.e. for all $t\geq t_0$ we have
	\begin{align}
	\left( \bm w_1(t) -\bm w_2(t)\right)^\top &\left(\left. \frac{\partial J_w(\bm x(t), \bm w)}{\partial \bm w}\right |_{\bm w_1(t)} \right. \nonumber \\
	&\quad\left.- \left. \frac{\partial J_w(\bm x(t), \bm w)}{\partial \bm w}\right |_{\bm w_2(t)}\right) > 0. \label{konvexitaetsforderung}
	\end{align}
	With $\bm w_1(t) = \bm w(t)$ and $\bm w_2(t) = \bm w^\star = \text{const.}$, we get
	\begin{align}
	\left. \frac{\partial J_w(\bm x(t), \bm w)}{\partial \bm w}\right |_{\bm w^\star} = \nullvec{}
	\end{align}
	and \eqref{konvexitaetsforderung} reads as 
	\begin{align}
		\left( \bm w(t) -\bm w^\star\right)^\top \left(\left. \frac{\partial J_w(\bm x(t), \bm w)}{\partial \bm w}\right |_{\bm w(t)}\right) > 0. \label{w-formel}
	\end{align}
Insertion of \eqref{adaptation} in \eqref{w-formel} yields
 \begin{align}
 \frac{1}{\alpha} 	\left( \bm w(t) -\bm w^\star\right)^\top \dot{\bm w}(t) < 0. \label{alpha-abstand}
 \end{align}
 With $\alpha >0$ and $\widetilde{\bm w}(t):=\bm w(t) - \bm w^\star$, \eqref{alpha-abstand} is equivalent to
 \begin{align}
 \left( \widetilde{\bm w}(t) \right)^\top \dot{\widetilde{\bm w}}(t) < 0. \label{kettenregel}
 \end{align}
 By using the chain rule, the left-hand side of \eqref{kettenregel} can be transformed to
 \begin{align}
 \frac 12 \cdot \frac{\mathrm d}{\mathrm dt} \left\{ (\widetilde{\bm w}(t))^\top \widetilde{\bm w}(t) \right\} < 0, \label{fast-am-ziel}
 \end{align}
 where $(\widetilde{\bm w}(t))^\top \widetilde{\bm w}(t) = \| \widetilde{\bm w}(t) \|_2^2$. Multiplying \eqref{fast-am-ziel} by two and applying the square root on both sides, we finally obtain
 \begin{align}
\frac{\mathrm d}{\mathrm dt}  \| \widetilde{\bm w}(t) \|_2 < 0,
 \end{align}
 i.e.  the distance $\| {\bm w}(t) - \bm w^\star \|_2$ strictly monotonically decreases with time for all $t \geq t_0$. This results in
 \begin{align}
 \lim_{t \to \infty} \left \| \bm w(t) - \bm w^\star \right\|_2 = 0.
 \end{align}
\qed \end{pf}
\begin{prop}\label{prop:irgendwann-in-W-plus}
		Let $V^\star(\bm x) = H(\bm x) + (\bm w^{\star})^\top \bm \Phi(\bm x)$ be the value function of optimization problem \eqref{original-OP} and let the conditions of Proposition~\ref{prop:konvexitat} hold with $\bm w_0 \in \mathcal N$. Then there exists a $T\geq t_0$ such that $\forall t > T \, : \,\bm w(T) \in \mathcal W^+$, 
		i.e. $\bm w(t)$ will remain in $\mathcal W^+$ after a finite amount of time and $\mathcal W^+$ is a positive invariant set for $t > T$.
		\end{prop}
	\begin{pf}
			With Lemma \ref{wastDrinne} and Lemma \ref{lemma:wPlusOpen}, $\bm w^\star \in \mathrm{int}\, \mathcal W^+$.
		This means that there exists an $\varepsilon > 0$ such that the ball $\mathcal B(\bm w^\star, \varepsilon) = \left\{\bm w \in \mathds R^r \, : \, \| \bm w - \bm w^\star \|_2 \leq \varepsilon\right\}$ lies completely within $\mathcal W^+$:
		\begin{align}
		\exists \, \varepsilon > 0 : && \mathcal B(\bm w^\star, \varepsilon) \subseteq \mathcal W^+
		\end{align}
		According to Proposition~\ref{streng-monoton-fallend}, $\left \| \bm w(t) - \bm w^\star \right\|_2$ is strictly decreasing with time. Consequently, there is a $T \geq t_0$ such that $\| \bm w(T) - \bm w^\star \|_2 = \varepsilon$,
		i.e. $\bm x,$ intersects the surface of the ball. Since $\left \| \bm w(t) - \bm w^\star \right\|_2$ is strictly decreasing, $\bm w(t)$ will remain within the ball and thus within $\mathcal W^+$ for all $t>T$.
	\qed \end{pf}
With this in mind, we can prove that $V(\bm x(t),\bm w(t))$ is indeed a positive-definite function:
\begin{prop}\label{prop:Vposdefinit}
		Let $V^\star(\bm x) = H(\bm x) + (\bm w^{\star})^\top \bm \Phi(\bm x)$ be the value function of optimization problem \eqref{original-OP} and let the conditions of Proposition~\ref{prop:konvexitat} hold with $\bm w_0 \in \mathcal N$. Then
		\begin{align}
		\forall \; t \geq t_0: && V(\bm x(t), \bm w(t)) \succ 0.
		\end{align}
		\end{prop}
		\begin{pf}
			According to Proposition~\ref{prop:irgendwann-in-W-plus}, there exists a $T \geq t_0$ such that 
			\begin{align}
			\forall t > T \, :  && \bm w(T) \in \mathcal W^+. \label{eq:remain-W}
			\end{align} Since $\mathcal W^+$ is the set of parameters $\bm w$ where $V(\bm x, \bm w)$ is positive-definite for $\emph{all}$ $\bm x \in \mathds R^n$, \eqref{eq:remain-W}
			implies that
			\begin{align}
			\forall \; t > T: && V(\bm x(t), \bm w(t)) \succ 0.
			\end{align}
	With $V(\bm x(t_0), \bm w(t_0))= V(\bm x_0, \nullvec{r})= H(\bm x_0) \succ 0$ and due to the fact that $V(\bm x(t), \bm w(t)$ is continuous and $\dot V(\bm x(t), \bm w(t))$ is monotonically decreasing according to Proposition~\ref{dot-V-decreasing}, $V(\bm x(t), \bm w(t) \succ 0$ for all $t\geq t_0$.
	\qed \end{pf}
As a consequence of Propositions \ref{dot-V-decreasing} and \ref{prop:Vposdefinit}, $V(\bm x(t),\bm w(t))$ is a suitable Lyapunov function. With this, we are ready to formulate the main statement of this paper regarding stability and asymptotic stability of the closed-loop equilibrium:
\begin{thm}
		Let $V^\star(\bm x) = H(\bm x) + (\bm w^{\star})^\top \bm \Phi(\bm x)$ be the value function of optimization problem \eqref{original-OP} and let the conditions of Proposition~\ref{prop:konvexitat} hold with $\bm w_0 \in \mathcal N$.
		Then $\bm x=\nullvec{}, \bm w = \bm w^\star$ is a stable equilibrium of \eqref{closed-loop}.
		If additionally one of the following conditions holds
	\begin{enumerate}
		\item The autonomous system $\dot{\bm x} = (\bm J-\bm R)\frac{\partial H}{\partial \bm x}$ is asymptotically stable with respect to the origin $\bm x = \nullvec{}$,
		\item $\bm G$ has full rank,
	\end{enumerate}
	then $\bm x=\nullvec{}, \bm w = \bm w^\star$ is an asymptotically stable equilibrium of \eqref{closed-loop}.
\end{thm}
\begin{pf}
	According to Proposition~\ref{prop:Vposdefinit}, $V(\bm x, \bm w)$ is positive-definite and according to Proposition~\ref{dot-V-decreasing}, $\dot V(\bm x, \bm w)$ is negative-semidefinite.
	As such, $V(\bm x, \bm w)$ is a Lyapunov function for the equilibrium $(\nullvec{}, \bm w^\star)$ of \eqref{closed-loop}, which is consequently a stable equilibrium.
	
	To prove asymptotic stability of $(\nullvec{} ,\bm w^\star)$, recall Proposition \ref{streng-monoton-fallend}, which states that $\bm w(t)$ converges strictly monotonically to $\bm w^\star$.
	Now let $\mathcal X^0 = \{ \bm x \in \mathds R^n: \dot V (\bm x(t),\bm w^\star)=0\}$ be the set of states where $V(\bm x, \bm w)$ is constant and $\bm w=\bm w^\star$. With regard to the individual summands in \eqref{dot_v_result}, we get
	\begin{align}
	\mathcal X^0 = \left\{ \bm x \in \mathds R^n:\right.& (\fl' = 0) \wedge ((\Ql' = 0) \vee (\Sl' = 0)),  \nonumber \\
	&\left.\bm w = \bm w^\star\right\}. \label{ix-null}
	\end{align}
	With $Q_\Upsilon' = r(\bm x)\succ 0$, the condition $\Ql' = 0$ is equivalent to $\bm x = \nullvec{}$, which implies $\fl'=0$. Accordingly, we can simplify \eqref{ix-null} to
	\begin{align}
		\mathcal X^0 = \left\{ \bm x \in \mathds R^n:\right.& (\bm x = \nullvec{}) \vee ((\fl' = 0) \wedge (\Sl' = 0)),  \nonumber \\
	&\left.\bm w = \bm w^\star\right\}. \label{ix-null-2}
	\end{align}
	From \eqref{ix-null-2}, conditions (1) and (2) of the Theorem are then obtained as follows:
	\begin{enumerate}
		\item According to LaSalle's invariance principle, all trajectories $\bm x(t)$ with $\dot V(\bm x(t),\bm w(t))=0$ converge to the largest invariant set contained in $\mathcal X^0$.
		Since $\bm S(\bm x)$ is positive-definite, $\Sl' = 0$ implies $\bm G^\top(\bm x) \frac{\partial V(\bm x, \bm w)}{\partial \bm x}= \nullvec{}$ (see \eqref{FormelFuerSUpsilon}). Bearing in mind \eqref{eq:closed_loop_eingang}, this leads to $\bm u= \nullvec{}$. Due to the assumption that the autonomous system is asymptotically stable with respect to $\bm x = \nullvec{}$,  the largest invariant set in $\mathcal X^0$ is a point. Thus $\bm x=\nullvec{}, \bm w = \bm w^\star$ is an asymptotically stable equilibrium of \eqref{closed-loop}.					
		\item if $\bm G(\bm x)$ has full rank, then $\bm G(\bm x) \bm S^{-1}(\bm x) \bm G^\top(\bm x) \succ 0$ and hence $\Sl' = 0$ only holds for $\bm x = \nullvec{}$. Thus $\mathcal X^0 = \{ \nullvec{}\}$, which implies that $\bm x=\nullvec{}, \bm w = \bm w^\star$ is an asymptotically stable equilibrium of \eqref{closed-loop}.	\qed
	\end{enumerate}
\end{pf}
\begin{rem}
To improve the convergence speed of $\widetilde{\bm w}(t)$, the gradient descent \eqref{adaptation} may be replaced by the continuous version of Newton's method \cite{Milano}:
\begin{align}
\dot{\bm w} = - \alpha \cdot \left(\frac{\partial^2 J_w}{\partial \bm w^2}\right)^{-1} \frac{\partial J_w}{\partial \bm x}, && \alpha > 0.
\end{align}
If ill-conditioning of the Hessian does not allow the numerical calculation of the inverse, there is a broad literature on alternative formulations of the Newton descent direction, such as Regularized Newton's Method \cite{Polyak.2009} or the pseudo-inverse formulation
\begin{align}
\dot{\bm w} = - \alpha \cdot \left(\frac{\partial^2 J_w}{\partial \bm w^2}\right)^{+} \frac{\partial J_w}{\partial \bm x}, && \alpha > 0.
\end{align}
Besides the ``classic'' least-squares pseudoinverse \cite{Golub.1965}, there exists a large number of advanced approaches based on singular value decomposition, e.g. truncated pseudoinverse or damped least-squares pseudoinverse, see \cite{Mullins.1995} for a discussion on alternative formulations.
\end{rem}
\section{Example}
In the following section, we apply the presented method to a linear and a nonlinear optimization problem of form \eqref{original-OP}, for which the value functions are explicitly known. Moreover, we compare the performance of our method with the performance of the ``exact'' optimal controller resulting from the value function.
\subsection{Linear Example}
At first, we consider the optimization problem \eqref{original-OP} with linear ISO-PHS dynamics
\begin{subequations}\label{linear-op}
\begin{IEEEeqnarray}{C/l}
\min_u & \frac{1}{2} \int_{\tini}^{\infty} \big(\bm x^\top \begin{bmatrix} 100 & 0 \\ 0 & 1\end{bmatrix} \bm x + u^2 \big) \dt	\\
\text{s.t.} &
\dot{\bm x} = \left[ \begin{bmatrix} 0 & -1 \\ 1 & 0 \end{bmatrix} - \begin{bmatrix} 1 & 0 \\ 0 & 1 \end{bmatrix} \right] \frac{\partial H(\bm x)}{\partial \bm x} + \begin{bmatrix} 1 \\ 0 \end{bmatrix} u +d,
\IEEEeqnarraynumspace \label{example-sys-linear}
\end{IEEEeqnarray}
\end{subequations}
where $H(\bm x) = \frac 12 x_1^2 + \frac 12 x_2^2$. Since $\rank (\Rphs) = 2$, the condition of Corollary \ref{cor:linearPHSCLF} is fulfilled and $H(\bm x)$ is a CLF. 

The exact solution for the value function can be calculated a priori to
\begin{align}
V^\star(\bm x) = \frac 12\bm x^\top \underbrace{\begin{bmatrix} \num{8.97697} & -\num{0.730021} \\ -\num{0.730021} & \num{0.963556}\end{bmatrix}}_{\bm P} \bm x, \label{Riccati}
\end{align}
where $\bm P$ is the Riccati matrix associated to \eqref{linear-op}.
The basis functions are chosen to $\bm \Phi(\bm x) = \begin{bmatrix} x_1^2 & x_1x_2  & x_2^2\end{bmatrix}^\top$.
By comparing \eqref{ansatz} and \eqref{Riccati}, we get
\begin{align}
\bm w^\star = \begin{bmatrix} \num{3.988485} & -\num{0.730021} & -\num{0.018222} \end{bmatrix}^\top.
\end{align}
The system is initialized at $\bm x_0 = \begin{bmatrix} 1 & 1 \end{bmatrix}^\top$, $\bm w_0 = \begin{bmatrix} 0 & 0 \end{bmatrix}^\top$. The shifts in $J_w$ (see \eqref{eq:linearkombinationNeu}) are set to $\bm c_1=\begin{bmatrix} 0 & 0 \end{bmatrix}^\top$, $\bm c_2=\begin{bmatrix} 1 & 0 \end{bmatrix}^\top$, $\bm c_3=\begin{bmatrix} 0 & 1 \end{bmatrix}^\top$, $\bm c_4=\begin{bmatrix} 1 & -1 \end{bmatrix}^\top$, and $\alpha$ is set to $\num{0,01}$.
\begin{figure}
	\centering
	\includegraphics[width=\columnwidth]{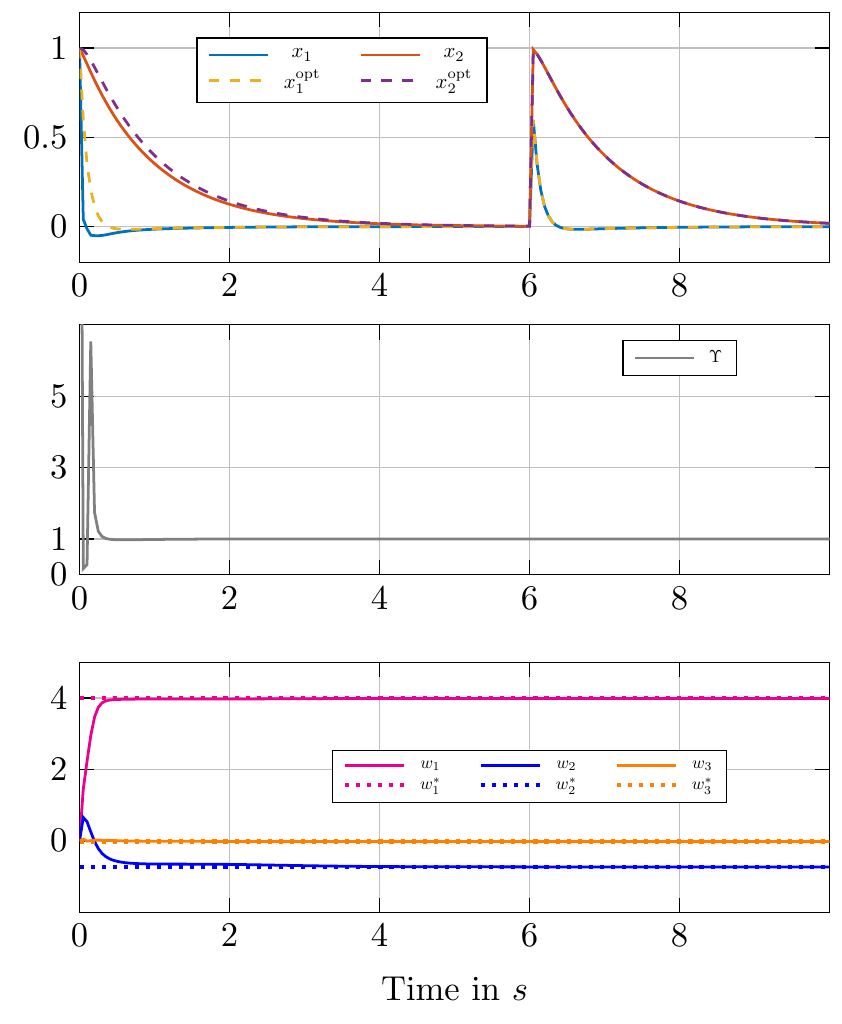}
	\caption{Adaptive (solid) and optimal controller (dashed) for optimization problem \eqref{original-OP} with linear ISO-PHS dynamics.}
		\label{fig:ex-linear}
\end{figure}

Fig. \ref{fig:ex-linear} shows the trajectories of $\bm x(t)$, $\Upsilon(t)$ and $\bm w(t)$, with the dashed curves indicating the optimum values associated to the Riccati solution.
It can be seen that after \SI{0,5}{\second}, both $\Upsilon$ and $\bm w$ have reached their optimal values. However, due to the learning process of $\w$, the trajectory of $\bm x$ does not match the Riccati solution $\bm x_{\mathrm R}$. After $t = \SI{6}{\second}$, when the adaption is finished, an additive disturbance input $d = \delta(t-6)$ is applied in order to evaluate the learning process. After the disturbance, $\Upsilon$ remains at $1$, $\bm w$ remains at $\bm w^\star$, and the closed-loop trajectory of $\bm x$ is identical to the Riccati solution $\bm x_{\mathrm R}$. Thus, it can be seen that the proposed controller converges to the optimal solution once the adaptation process of the value function parameters is finished. Even after the disturbance, the parameters $\bm w$ remain at their optimum value $\bm w^\star$ due to the strict convexity of the function in \eqref{adaptation} used for adaption (c.f. Proposition~\ref{prop:konvexitat}).
\subsection{Nonlinear Example}
Next, we consider the following nonlinear optimization problem
\begin{subequations}\label{nl-op}
\begin{IEEEeqnarray}{C.l}
\min_u & \frac{1}{2} \int_{\tini}^{\infty} \big(\bm x^\top \begin{bmatrix} 8+8x_1+16x_2 & 0 \\ 0 & 8\end{bmatrix} \bm x + u^2 \big) \dt	\label{nl-op1}\\
\text{s.t.} &\dot{\bm x} = \left[ \begin{bmatrix} 0 & 3 \\ -3 & 0 \end{bmatrix} - \begin{bmatrix} 1+x_2^2 & 1 \\ 1 & 2 \end{bmatrix} \right] \frac{\partial H(\bm x)}{\partial \bm x} + \begin{bmatrix} x_2 \\ 0 \end{bmatrix} u 
\IEEEeqnarraynumspace
\label{nl-op2}
\end{IEEEeqnarray}
\end{subequations}
with Hamiltonian $H(\bm x) = \frac 12 x_1^2 + \frac 12 x_2^2$. 
Following the lines of \cite{nevistic1996constrained}, the value function for this specific optimization problem is known to be 
\begin{align}\label{eq:valuef}
V^\star(\bm x) = 2x_1^2 + x_2^2.
\end{align}
Since the dissipation matrix $\bm R$ is positive-definite for all $\bm x\in \mathds R^2$, the condition of Proposition \ref{prop:HistCLF} is satisfied and $H(\bm x)$ is a CLF. 


With the choice
$
\bm \Phi (\bm x) = \begin{bmatrix} x_1^2  & x_1x_2 & x_2^2 \end{bmatrix}^\top
$,
and by comparing \eqref{ansatz} and \eqref{eq:valuef}, the optimal weighting factors are
$
\bm w^\star = \begin{bmatrix} \num{1,5} & \num{0} & \num{0,5} \end{bmatrix}^\top.
$
Again, the system is initialized at $\bm x_0 = \begin{bmatrix} 1 & 1 \end{bmatrix}^\top$, $\bm w_0 = \begin{bmatrix} 0 & 0 \end{bmatrix}^\top$ and the disturbance input $d = \delta(t-6)$ is added. 
The shifts in $J_w$ are set to $\bm c_1=\begin{bmatrix} 0 & 0 \end{bmatrix}^\top$, $\bm c_2=\begin{bmatrix} -1 & 0 \end{bmatrix}^\top$, $\bm c_3=\begin{bmatrix} 0 & -1 \end{bmatrix}^\top$, $\bm c_4=\begin{bmatrix} 1 & -1 \end{bmatrix}^\top$, and $\alpha$ is set to \num{0,02}.
\begin{figure}
	\centering
	\includegraphics[width=\columnwidth]{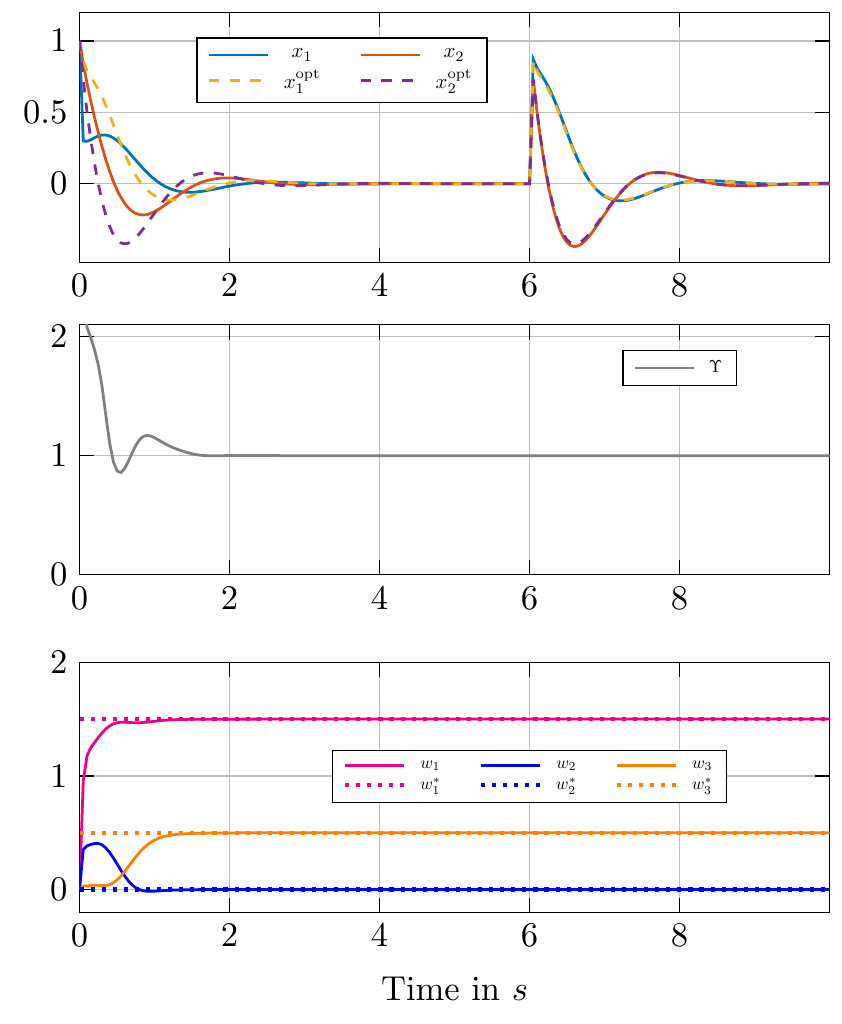}
	\caption{Adaptive (solid) and optimal controller (dashed) for optimization probleme \eqref{original-OP} with nonlinear ISO-PHS dynamics.}
		\label{fig:ex-nichtlinear}
\end{figure}

Fig. \ref{fig:ex-nichtlinear} shows the trajectories of $\bm x(t)$, $\Upsilon(t)$ and $\bm w(t)$, with the dashed curves indicating the optimum values associated to the optimal controller $\bm u = -\mm{S}^{-1}\bm G^\top \frac{\partial V^\star}{\partial x}$.
It can be seen that after $1s$, both $\Upsilon$ and $\bm w$ have reached their optimal values. Moreover, the trajectory of $\bm x$ converges to the optimal solution associated to the optimal controller and remains identical once the learning process is completed, even after the disturbance. 
Overall, these examples demonstrate that the proposed controller is capable of adapting the optimal controller parameters after a single learning phase.

\begin{figure}
	\centering
	\includegraphics[width=\columnwidth]{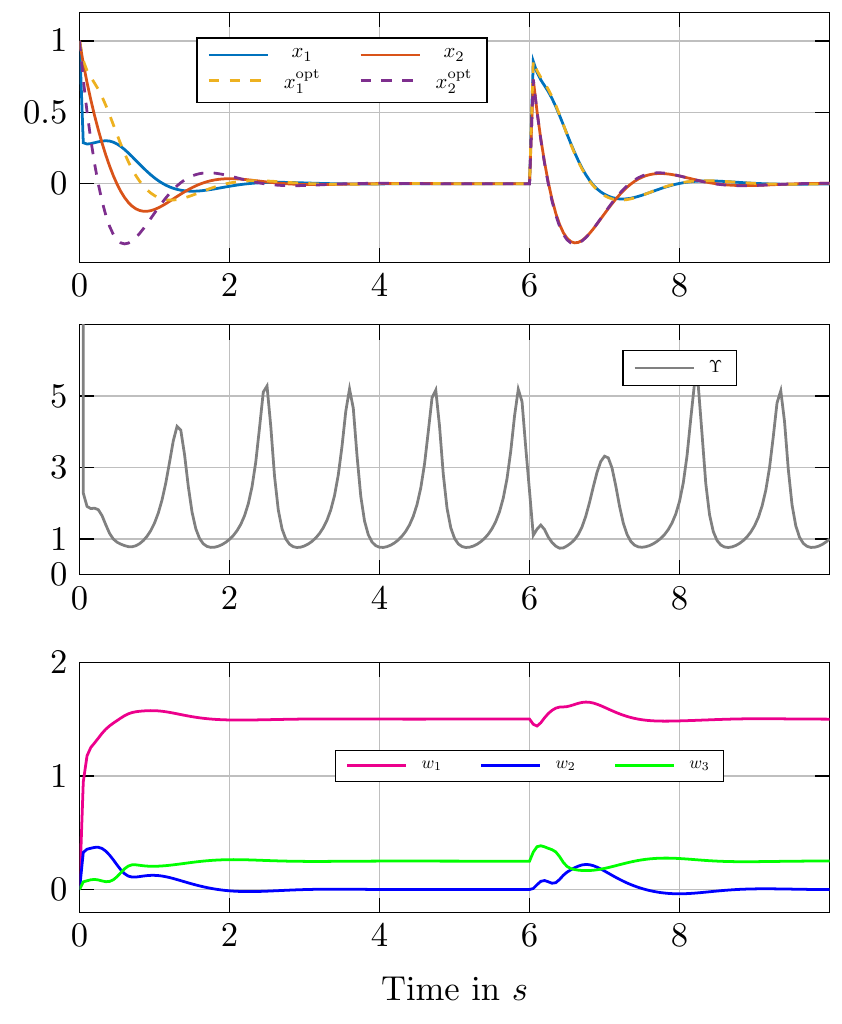}
	\caption{Adaptation if Assumption \ref{ass:optimal} is violated.}
		\label{fig:ex-nichtlinear-falsch}
\end{figure}

To investigate the effects of an incorrect choice of basis functions, we repeat the simulation for optimization problem \eqref{nl-op} using
$
\bm \Phi ' (\bm x) = \begin{bmatrix} x_1^2 & x_1x_2 & x_2^4 \end{bmatrix}^\top,
$
i.e. $H(\bm x)+\bm w^\top \bm \Phi '(\bm x)$ does not fit the structure of $V^\star(\bm x)$ and hence Assumption \ref{ass:optimal} is violated. The results are shown in Fig. \ref{fig:ex-nichtlinear-falsch}.
The trajectory of $\Upsilon$ shows a remarkable and distinct oscillatory behaviour. However, after about $3s$, the weighting factors $\bm w$ converge to a certain value $\bm w^\diamond$.
After the disturbance at $t=\SI{6}{\second}$, the weighting factors do not remain at their previous values. 

Both the oscillation of $\Upsilon (t)$ and the fluctuation of $\bm w(t)$ after the disturbance imply that $H(\bm x) + (\bm w^\diamond)^\top\bm \Phi'(\bm x)$ is not equal to the value function (cf. Remark~\ref{rem:suboptimal}). 
However, the results show that even if $\bm \Phi '(\bm x)$ is not accurate, the proposed controller is able to learn suitable weighting factors for a suboptimal control. Furthermore, the oscillation of $\Upsilon(t)$ can be interpreted as an indicator of suboptimality for the chosen set of functions in $\bm \Phi ' (\bm x)$, as discussed in Remark \ref{rem:suboptimal}.
\section{Conclusion}
In this paper, we have introduced a time-continuous adaptive feedback controller 
for dynamic optimization with generalized Lagrangian performance indices and general time-continuous ISO-PHSs.
In particular, we stated necessary and sufficient conditions under which the Hamiltonian $H(\bm x)$ is a CLF.
As a consequence, the initial value function guess $V(\bm x)=H(\bm x)$ allows to deploy an admissible controller which is already stabilizing.
Based on this initial guess, we proposed a gradient-based continuous learning procedure for the extended CLF $V(\bm x, \bm w) = H(\bm x)+ \bm w^\top \bm \Phi(\bm x)$ with the aim of approximating the value function $V^\star(\bm x)$. 
We proved (asymptotic) stability of the closed-loop system equilibrium $(\nullvec{}, \bm w^\star)$.
Finally, we investigated our theoretical findings by means of a linear and a nonlinear simulation example.


Although a reasonable choice of basis functions is nontrivial, simulations show that even if $\bm \Phi(\bm x)$ is inaccurate, the controller is able to stabilize the system, providing near-optimal solution trajectories. Furthermore, the optimality of the computed control law can be assessed for the case of a bad choice of basis functions via $\Upsilon(t)$. 
However, a rigorous perturbation analysis for systems where parameterization $\Phi(\bm x)$ does not fit the structure of $V^\star(\bm x)$ remains an open research question. 
\appendix
\section{Appendix}
\subsection{Lemma \ref{lemma:eigenwert_und_kern}}\label{app:Kernlemma}
\begin{lem}\label{lemma:eigenwert_und_kern}
	Consider a symmetric, positive-semidefinite matrix $\mm{M} \in \R^{n\times n}$. Then
	\begin{align}
	\x^\top \mm{M} \x = 0 && \Longleftrightarrow && \bm x \in \ker\{\mm{M}\}. \label{quadrat-null}
	\end{align}
\end{lem}
\begin{pf}
	``$\Longleftarrow$'': trivial. 
	
	``$\Longrightarrow$'':
	For a positive-semidefinite, symmetric matrix $\mm{M}$ with real entries all eigenvectors are orthogonal \cite[Th.~7.2.1]{anton2013}. Consequently, each vector $\x$ can be expressed as a linear combination of the eigenvectors $\mv{v}_i$, 
	\begin{align}\label{eq:ew}
	\x = \sum_{i=1}^{n}\varsigma_i \mv{v}_i
	\end{align}
	with $\varsigma_i \in \mathds R$. It follows that the product of $\mm{M}$ and $\x$ can be written as
	\begin{align}
	\mm{M}\x = \sum_{i=1}^{n}\varsigma_i \mm{M}\mv{v}_i = \sum_{i=1}^{n}\varsigma_i \lambda_i \mv{v}_i . \label{that}
	\end{align} 
	With \eqref{that} and taking into account that $\mv{v}_i^\top \mv{v}_j = 0$, $i \neq j$ due to the orthogonality, \eqref{quadrat-null} can be written as
	\begin{align}
	\bigg( \sum_{i=1}^{n}\varsigma_i \mv{v}_i \bigg)^\top \sum_{i=1}^{n}\varsigma_i \lambda_i \mv{v}_i = \sum_{i=1}^{n} \varsigma_i \lambda_i \mv{v}_i^\top \mv{v}_i = 0,
	\end{align}
	which is not equal to zero unless all the eigenvalues $\lambda_i$ of the eigenvectors used for $\x = \sum_i \varsigma_i\mv{v}_i$ are zero. Hence, all the solutions to equation $\x^\top \mm{M} \x = 0$ are spanned by eigenvectors corresponding to an eigenvalue zero. The vector space $\ker\{\mm{M}\} = \{\x \mid \mm{M}\x = \nullvec{n}\}$ with \eqref{eq:ew} can be written as 
	\begin{align}
	\mm{M}\x = \sum_{i=1}^{n}\varsigma_i\mm{M}\mv{v}_k = \sum_{i=1}^{n}\varsigma_i\lambda_k\mv{v}_k = 0,
	\end{align}
	which is also spanned by the eigenvectors corresponding to eigenvalue zero. Hence, it is easy to see that all the eigenvectors corresponding to the eigenvalue zero constitute exactly the vector space $\ker\{\mm{M}\}$ and also the vector space constituted by $\x^\top \mm{M} \x = 0$. 
	\qed \end{pf}
\subsection{Lemma \ref{lemma:wPlusOpen}}\label{ch:openLemma}
\begin{lem}\label{lemma:wPlusOpen}
	$\mathcal W^+$ is an open set, i.e. for each $\bm w \in \mathcal W^+$ there is an $\varepsilon >0$ such that each $\bm w' \in \mathds R^r$ with $\| \bm w - \bm w' || < \varepsilon$ lies within $\mathcal W^+$.
\end{lem}
\begin{pf}
	To prove that $\mathcal W^+$ is an open set, two auxiliary sets $\mathcal Y^+$ and $\mathcal Z$ are introduced which can be proved to be open more easily. Finally, openness of $\mathcal W^+$ is concluded by a canonical projection of $\mathcal Z$:

Let $\mathcal Y^+ \subseteq \mathds R^n \times \mathds R^r$ be the set of all $(\bm x,\bm w)$ with $V(\bm x, \bm w) >0$:
\begin{align}
\mathcal Y^+ = \left\{(\bm x, \bm w) \in \mathds R^n \times \mathds R^r \; : \; V(\bm x, \bm w) > 0\right\}
\end{align} 
Thus $\mathcal Y^+$ is the preimage of $\mathds R^+$, i.e. $\mathcal Y^+=V^{-1}(\mathds R^+)$. Each preimage of a continuous function is open whenever the corresponding image is open \cite[Theorem 2.9]{Crossley.2005}.
Since $V(\bm x,\bm w)$ is continuous and $\mathds R^+$ is open, $\mathcal Y^+$ must also be open. Furthermore, we note that $\mathcal Y^+$ is nonempty due to the fact that $V(\bm x, \nullvec{r})=H(\bm x) > 0$ is fulfilled by definition for all $\bm x \in \mathds R^n\backslash \{\nullvec{}\}$.

Now consider the \emph{$\mathcal Y^+$-inner cylinder} (see Fig. \ref{fig:cylinder}) with
\begin{align}
\mathcal Z(\underline{\bm w}, \overline{\bm w}):= &\left\{(\bm x, \bm w) \in \mathcal Y^+ \; : \right. \nonumber \\
&\left. \bm w \in ]\underline{\bm w}, \overline{\bm w}[  , (\widetilde{\bm x}, \bm w) \in \mathcal Y^+ \; \forall \; \widetilde{\bm x} \in \mathds R^n \backslash \{\nullvec{}\}  \right\}.
\end{align}
Obviously $\mathcal Z(\underline{\bm w}, \overline{\bm w}) \subseteq \mathcal Y^+$ and according to Lemma \ref{lemma:zylinder} (see below), $\mathcal Z(\underline{\bm w}, \overline{\bm w})$ is open for each pair $(\underline{\bm w}, \overline{\bm w}) \in \mathds R^r \times \mathds R^r$.

With the help of $\mathcal Z(\underline{\bm w}, \overline{\bm w})$, we can define the \emph{maximum $\mathcal Y^+$-inner cylinder}
\begin{align}
\mathcal Z^+ = \bigcup_{\underline{\bm w}, \overline{\bm w} \in \mathds R^r} \mathcal Z (\underline{\bm w}, \overline{\bm w}).
\end{align}
as the union of all possible cylinders $\mathcal Z (\underline{\bm w}, \overline{\bm w})$, see Fig. \ref{fig:cylinder}. 
\begin{figure}
	\centering
	\includegraphics[width=0.4\textwidth]{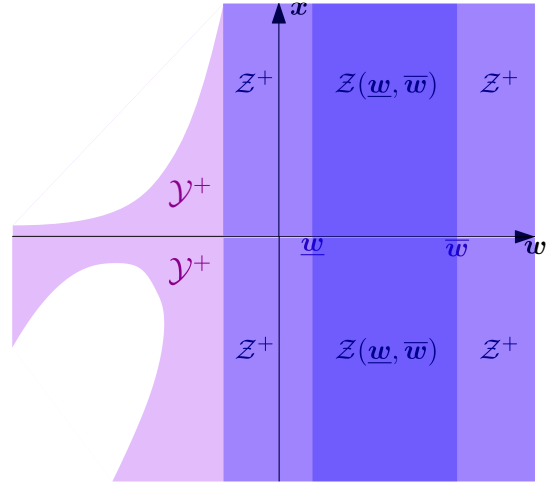}
	\caption{Ilustration of $\mathcal Y^+$ with an inner cylinder $\mathcal Z(\underline{\bm w}, \overline{\bm w})\subseteq \mathcal Y^+$.}
	\label{fig:cylinder}
\end{figure}
Since the union of open sets is open \cite[Theorem 1.1.9]{Singh.2019}, also $\mathcal Z^+$ is open.
With the canonical projection 
\begin{align}
\mathrm{proj}: \mathds R^n \times \mathds R^r \to \mathds R^r,
\end{align}
the set $\mathcal W^+$ can be interpreted as
\begin{align}
\mathcal W^+ = \mathrm{proj}(\mathcal Z^+),
\end{align}
i.e. the canonical projection of $\mathcal Z^+$ in $\mathds R^r$. Since projection maps are open maps \cite[p.~5]{Deo.2018}, $\mathcal W^+$ is an open set.
\qed \end{pf}
\begin{lem} \label{lemma:zylinder}
	The $\mathcal Y^+$-inner cylinder
	\begin{align}	\mathcal Z(\underline{\bm w}, \overline{\bm w}):= &\left\{(\bm x, \bm w) \in \mathcal Y^+ : \right.\nonumber \\
	&\quad\left. \bm w \in ]\underline{\bm w}, \overline{\bm w}[  , (\widetilde{\bm x}, \bm w) \in \mathcal Y^+ \; \forall \; \widetilde{\bm x} \in \mathds R^n \backslash \{\nullvec{}\}  \right\} \label{zylinder}
	\end{align}
	with $\underline{\bm w}, \overline{\bm w} \in \mathds R^r$ is an open set.
\end{lem}
\begin{pf}
	If $\underline{\bm w} \geq \overline{\bm w}$, then the interval $\left]\underline{\bm w}, \overline{\bm w}\right[$ is improper, thus $\mathcal Z(\underline{\bm w}, \overline{\bm w})=\emptyset$. Since empty sets are trivially open, the proof is complete. 
	
	For $\underline{\bm w} <  \overline{\bm w}$ let $(\bm x', \bm w') \in \mathcal Z (\underline{\bm w}, \overline{\bm w})$ be an arbitrary point within the cylinder (see Fig. \ref{fig:cylinderBall}). Now let
	\begin{align}
	\bm x'' &= \text{col}_i\{| x_i' \| \}, && i = 1,\ldots, n\\	
	\bm w'' &= \text{col}_j\{\min \{ \overline w_j - w_j', w_j' - \underline w_j\} \}, && j = 1, \ldots, r,
	\end{align}
	where $\bm x''$ denotes the componentwise distances between $x_i'$ and $0$ and  $\bm w''$ denotes the componentwise distances between $w_j'$ and the lower or upper bounds $\underline w_j$ or $\overline w_j$, respectively. 
	\begin{figure}
		\centering
		\includegraphics[width=0.4\textwidth]{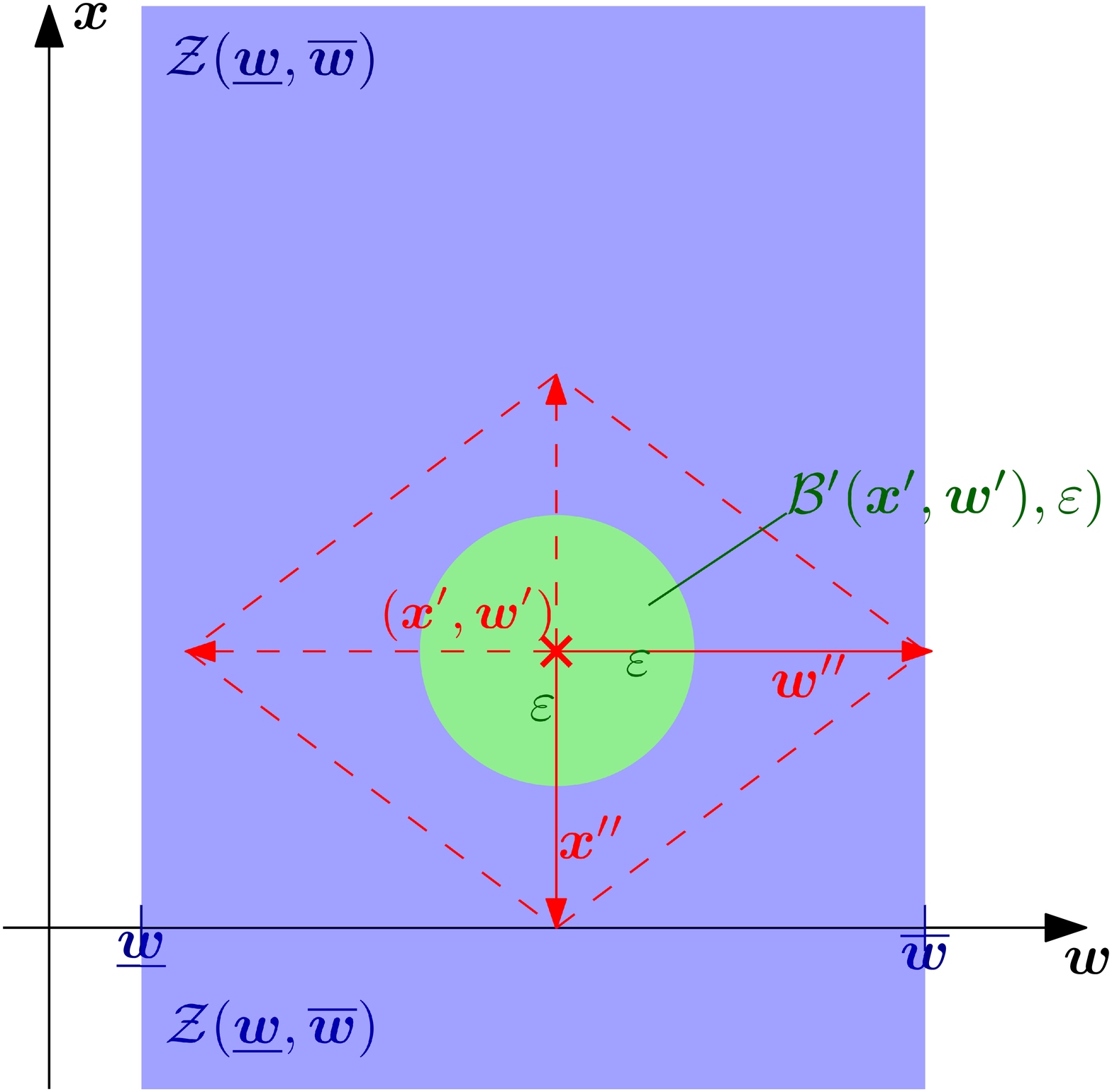}
		\caption{Each inner cylinder $\mathcal Z(\underline{\bm w}, \overline{\bm w})$ is an open set.}
		\label{fig:cylinderBall}
	\end{figure}
	From Fig. \ref{fig:cylinderBall} it can be seen that by definition of $\bm x''$ and $\bm w''$:
	\begin{align}
	(\bm x' + \frac 12 \bm x'', \bm w') &\in \mathrm{int}\,\mathcal Z(\underline{\bm w}, \overline{\bm w}), \\
	(\bm x' - \frac 12 \bm x'', \bm w') &\in \mathrm{int}\,\mathcal Z(\underline{\bm w}, \overline{\bm w}), \\
	(\bm x' \bm w'+ \frac 12 \bm w'') &\in \mathrm{int}\,\mathcal Z(\underline{\bm w}, \overline{\bm w}), \\
	(\bm x', \bm w' - \frac 12 \bm w'') &\in \mathrm{int}\,\mathcal Z(\underline{\bm w}, \overline{\bm w}).
	\end{align}
Thus, it is obvious that we can always construct an open ball
	$
	\mathcal B'((\bm x', \bm w'), \varepsilon)
	$
	around $(\bm x', \bm w')$ with radius
	\begin{align}
	\varepsilon = \frac 12 \cdot \left\| \begin{bmatrix} \bm x'' \\ \bm w'' \end{bmatrix} \right\|_\infty
	\end{align}
	that lies completely in $\mathcal Z(\underline{\bm w}, \overline{\bm w})$. Hence $\mathcal Z(\underline{\bm w}, \overline{\bm w})$ is an open set.
\qed \end{pf}
\bibliography{root}
\end{document}